# A proof that $h_1, h_2 \leq h_0$ for any h-basis $A_3$

MICHAEL CHALLIS

*Storey's Cottage, 3 Church Lane, Whittlesford, Cambridge CB2 4NX, UK*

26th January 1996

**History**

| | | |
|---|---|---|
| 0.01 | 21-Dec-95 | Document started. |
| 0.02 | 26-Jan-96 | First version complete; still relies on [1] for some proofs. |
| 0.03 | 13-Feb-96 | Proof of Lemma 10 and new proof of Lemma 12 added; this document is now self-contained (ie does not rely on [1]). |
| 0.04 | 16-Feb-96 | Start adding details of alternative approaches and additional information determined during the investigation; all identified by a smaller point size. |
| 0.05 | 18-Mar-96 | All complete except for the case $C_2=1$, $C_1>a_2/2$, $p=2$, $s>n'/2$, $s'\geq 3$, $k=2$ ... |
| 0.06 | 20-Mar-96 | Document completed |

**Abstract**

$A_k = \{1, a_2, \ldots a_k\}$ is an h-basis for X if every positive integer $\leq X$ can be expressed as the sum of no more than h values $a_i$; X(h) is called the h-range of the basis. $h_0$ is the smallest value of h for which $X(h) \geq a_k$, and $h_1$ is the smallest value for which $X(h+1) = X(h) + a_k$ for all $h \geq h_1$. $h_2 \geq h_1$ identifies a further "stabilisation" in the h-range - a definition is included in the body of the paper. It is known that $h_1, h_2 \leq h_0$ for h-bases $A_3$, but published proofs are complicated (see Ch. VIII of [4] for a discussion, where references [3] and [5] are given). This paper introduces the concept of a "stride generator" $A = \{1, a_2, a_3\}$ which, while sharing some of the properties of a basis $A_3$, is simpler to treat mathematically. We establish a relationship between stride generators and h-bases, and show that $h_1, h_2 \leq h_0$ follows immediately if the stride generator underlying a basis has a particular property - here called "canonicality". The proof is lengthy (with a number of special cases to consider), but the underlying principles remain simple.

**Contents**

# 1 Stride generators and h-bases

## 1.1 Introduction and definitions

Let $A = \{1, a_2, a_3\}$ be a set of integers with $1 < a_2 < a_3$; we write $a_3 = C_2 a_2 + C_1$ where $0 \leq C_1 < a_2$.

An *h-basis* $B(A, h)$ has the following properties:

We say x has an *h-representation* if $x = c_3 a_3 + c_2 a_2 + c_1$ for $c_i \geq 0$, $c_1 + c_2 + c_3 \leq h$.

The basis' *h-range* $X(h)$ is defined as one less than the smallest integer which has no h-representation.

We say that the basis is *admissible* if $X(h) \geq a_3$; the smallest value of h for which this is true is denoted $h_0$. In what follows, we consider only admissible bases.

It is easy to show that $X(h+1) \geq X(h) + a_3$ for all $h \geq h_0$, and that there is a value $h_1 \geq h_0 - 1$ beyond which equality obtains.

All values less than or equal to $X(h)$ have h-representations, and no value greater than $ha_3$ has one; there may or may not be a representation for a value $X(h) < x < ha_3$. It can be shown that there is a value $h_2 \geq h_1$ such that for all $h \geq h_2$:

x has no h-representation $\iff$ $(x+a_3)$ has no (h+1)-representation for all $X(h) < x < ha_3$

This paper proves the following for all admissible h-bases $B(A, h)$:

(1)  $X(h+1) = X(h) + a_3$

(2)  x has no h-representation $\iff$ $(x+a_3)$ has no (h+1)-representation for all $X(h) < x < ha_3$

In other words, $h_1 \leq h_0$ and $h_2 \leq h_0$.

A *stride-generator* $SG(A, n, p)$ has the following properties:

We say x has an *n-generation* if there exists $i \geq 0$ such that $x + ia_3 = c_2 a_2 + c_1$ for $c_i \geq 0$, $c_1 + c_2 \leq n+i$; such a generation is said to be of *order* i.

Every integer $0 \leq x < a_3$ has an n-generation of order $\leq p$.    (A)

At least one integer $0 \leq x < a_3$ has no n-generation of order $< p$.    (B)

At least one integer $0 \leq y < a_3$ has no (n-1)-generation of order $\leq p+1$.    (C)

We can think of a stride generator as a recipe for representing each value $ka_3 \leq x' < (k+1)a_3$ for sufficiently large k, since if x has an n-generation of order i then $x + ka_3$ has an (n+k)-representation provided that $k \geq i$. With this view, y is (one of) the most difficult values to generate, since $y + ka_3$ has an (n+k)-representation, but no (n+k-1)-representation - at least for $k \leq p+1$.

Any value y which satisfies condition (C) is called a *break* in the stride generator.

If there is no value j such that $y + ja_3 = c_2 a_2 + c_1$ is soluble for $c_2 + c_1 \leq (n-1) + j$, we say that y is a *canonical* break; otherwise, we say that y has *break order* q where $q > p + 1$ is the smallest value of j for which the above equation has a solution.

We say that a stride generator is *canonical* if all of its breaks are canonical.

Lemma 11 below clarifies the relationship between h-bases and stride generators; it shows that every h-base $B(A, h)$ with h-range X has an underlying stride generator $SG(A, h-k, p)$ with a break $y = X + 1 \pmod{a_3}$. It turns out that all underlying stride generators are canonical, and it is from this property that we deduce easily that $h_1, h_2 \leq h_0$.

Stride generators are best understood when represented as thread diagrams:

A *thread T(e, i) of order i* is a contiguous sequence of integers [c, d], d ≥ c, corresponding to a sequence of n-generations all of the same order i:

$c + ia_3 = ea_2$
$(c+1) + ia_3 = ea_2 + 1$
$\ldots$
$d + ia_3 = ea_2 + (d-c)$  where  $e + (d - c) = n + i$

We write:

str(T) = $c = ea_2 - ia_3$  - the *start* of the thread
end(T) = $d = (ea_2 - ia_3) + (n + i) - e$  - the *end* of the thread
len(T) = $(d - c) + 1 = (n + i) - e + 1$  - the *length* of the thread
ord(T) = $i$  - the order of the thread

A *thread diagram* is an (x, y) diagram in which every thread T(e, i) = [c, d] is represented by a horizontal line at height y = i running from x = c to x = d inclusive; this line is optionally labelled e. The diagram covers the range $0 \leq x < a_3$.

A value x is *covered* by a thread T if $c \leq x \leq d$.

If $T_1 = [c_1, d_1]$ and $T_2 = [c_2, d_2]$, then $T_1$ *covers* $T_2$ if $c_1 \leq c_2$ and $d_1 \geq d_2$.

A value x is *crossed* by a thread T if $c \leq x < d$; in other words, T crosses x if it covers both x and x+1.

Unless otherwise stated, we consider only threads which cover at least one value $0 \leq x < a_3$; in other words, threads which at least partly appear in the stride generator's thread diagram.

It is easy to see that the following is an equivalent definition of a stride generator in terms of its thread diagram:

Every value $0 \leq x < a_3$ is covered by some thread of order $i \leq p$.  (A)

At least one value $0 \leq x < a_3$ is not covered by any thread of order $i < p$.  (B)

At least one value $0 \leq y < a_3$ has the property that no thread of order $\leq p+1$ crosses y.  (C)

If there is no thread of any order that crosses y, the break is canonical; otherwise, its break order is that of the first thread to cross y.

Every set A has at least one stride generator, and sometimes several different ones; as an example, A = {1, 38, 97} has three stride generators SG(A, 19, 2), SG(A, 15, 4) and SG(A, 14, 6):

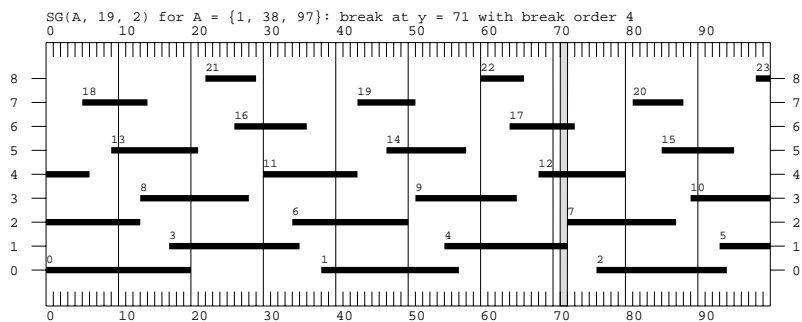

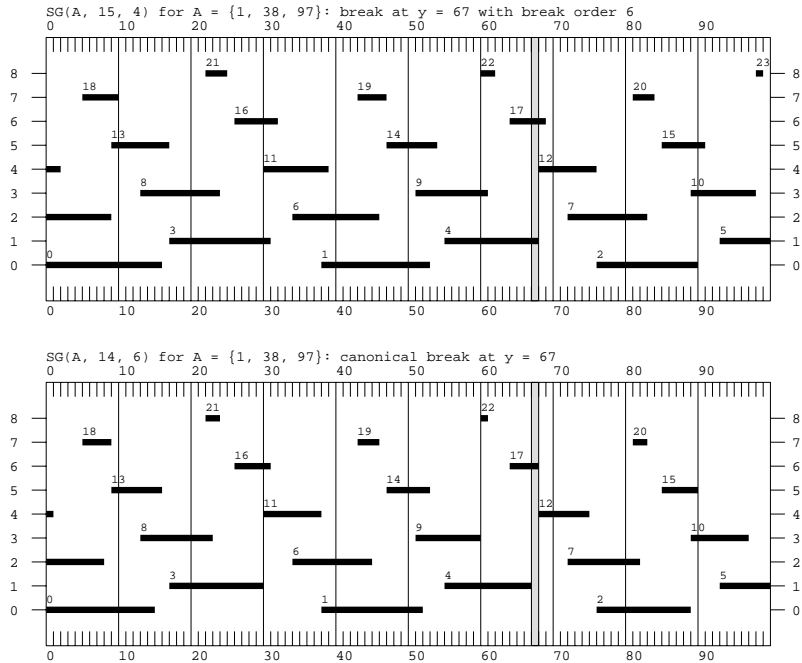

Some basic properties of threads (which can be seen in the diagrams above) are:

Threads of the same order recur at intervals of $a_2$; each one is one shorter than its predecessor.

More formally, if $T_1$ and $T_2$ are two consecutive threads of the same order, then $str(T_2) = str(T_1) + a_2$, and $len(T_2) = len(T_1) - 1$.

Threads whose orders differ by 1 are separated by $C_1$, and differ in length by $(C_2 - 1)$.

More formally, if $T_1 = T(e, i)$ and $T_2 = T(e+C_2, i+1)$ are two threads, then $str(T_1) = str(T_2) + C_1$, and $len(T_2) = len(T_1) - (C_2 - 1)$.

This means that any pattern of threads can be moved from one position in a thread diagram to another - either by moving the start position of each thread by a multiple of $a_2$, or altering the order of each thread by some constant - subject to every thread in the pattern retaining a positive length and order; I call this the *similarity property*.

More formally, suppose $T_1$ and $T_2$ are two threads related as follows:

$$ord(T_1) - ord(T_2) = x, \quad str(T_1) - str(T_2) = y, \quad len(T_1) - len(T_2) = z$$

where $T_1$ is at least as long as $T_2$ ($z \geq 0$). Let $T$ be any other thread with $len(T) > z$; then there exists a thread $U$ where:

$$ord(U) = ord(T) - x, \quad str(U) = str(T) - y, \quad len(U) = len(T) - z.$$

### 1.2 Properties of stride generators

Our first few lemmas prove some simple properties about the threads in the thread diagram for a stride generator $SG(A, n, p)$.

*Lemma 1*

$y \geq a_3 - a_2$ for any break $y$ in a stride generator.

Proof

Suppose the contrary, and consider $y' = y + a_2$; since $y' < a_3$, it must be covered by some thread $T = T(a, i)$. Now consider thread $T' = T(a-1, i)$. We have $str(T') = str(T) - a_2$ and $len(T') = len(T) + 1$; so $T'$ must cover both $y$ and $y+1$, and so $y$ cannot be a break.

*Page 5*

*Lemma 2*

There exists $x \geq a_3 - a_2$ that satisfies condition (B): that is, x is covered only by a thread of order p.

Proof

Let x satisfy condition (B), and suppose $x < a_3 - a_2$:

Consider $x' = x + a_2$. By condition (A), x' must be covered by some thread V' of order $j \leq p$, and hence x is covered by the thread V of order j with $str(V) = str(V') - a_2$; so $j = p$. Repeat until $x' \geq a_3 - a_2$.

*Lemma 3*

All threads of order $i \leq p$ exist.

Proof

If the last thread of order p exists, then all other threads of order p and all threads of lower order must also exist; so it is sufficient to show that there exists a thread T of order p that satisfies $a_3 - a_2 \leq str(T) < a_3$.

By Lemma 2 we may choose $x \geq a_3 - a_2$ that is covered by a thread T of order p. Suppose $str(T) < a_3 - a_2$; then $len(T) \geq 2$, and so there exists T' of order p with $str(T') = str(T) + a_2$ and so $str(T') \geq a_3 - a_2$.

*Lemma 4*

If there is no thread T(e, i) of order $i > p$, then any thread T(f, i) for $f < e$ is covered by some thread of order $j \leq p$.

(In other words, we can treat a *missing* thread as if it were a *covered* thread.)

Proof

Let $T' = T(e', i)$ for $e' < e$ be the first thread of order i - if any - that exists: $len(T') = 1$. Since T' is part of a stride generator, the value $x = str(T')$ must be covered by some thread $U' = T(g', j)$ for some $j \leq p$. Since $len(T') = 1$, this means that T' is covered by U' - and hence $T'' = T(e' - k, i)$ is covered by $U'' = T(g' - k, j)$ for any $k \geq 0$ as required.

*Lemma 5*

If a thread T of order i is covered by some other thread U of order $j < i$, then any thread V of order $k \geq i$ is covered by some thread V' of order $k' < i$.

Proof

By the similarity property, there is a thread $V_1$ of order $k_1 = k - (i-j)$ that covers V. If $k_1 < i$, the lemma is proved; otherwise we apply the similarity property repeatedly until we find thread $V' = V_n$ of order $k' = k_n = k - n(i-j)$ with $k' < i$ which covers V as required.

Lemma 4 and 5 together say that once we have found a thread of order i that does not exist or is covered by some other thread, we need only consider threads of order < i when looking for generations: for in such a case x has a generation if and only if it is covered by some thread of order $\leq i$.

*Lemma 6*

No thread of order $i \leq p$ is covered by any other thread.

Proof

If such a thread existed, every value $0 \leq x < a_3$ would be covered by some thread of order < i,

which is contrary to condition (B) for a stride generator.

One immediate corollory of Lemma 6 is that no two threads of order $\leq p$ can both start or both end in the same position in a stride generator.

*Lemma 7*

A stride generator is canonical if a thread of order $p+1$ does not exist, or is covered by some other thread.

Proof

By Lemma 4 and 5 this means that all threads of order $\geq p+1$ are covered by threads of order $< p$. Let y be any break in the stride generator, and suppose that it is crossed by a thread T of order $q > p+1$; then it must also be crossed by the thread U of order $< p$ that covers T, and so cannot be a break after all. So y is a canonical break, and the lemma is proved.

Later, we show the converse of the above Lemma: that for a canonical stride generator all threads of order $i > p$ are covered by threads of order $j \leq p$. From this we deduce that $h_2 = h_0$.

*Lemma 8*

Let $x \geq a_3 - a_2$ satisfy condition (B); then $x' = (x - C_1)$ is covered only by the thread $T(C_2-1, 0)$.

Proof

Suppose x' is covered by a thread $T(a, i)$ of order i. If $i > 0$, x will be covered by the thread $T(a-C_2, i-1)$ whose order $i-1 < p$; so i must be zero.

$a_3 - a_2 \leq x < a_3 \Rightarrow (C_2-1)a_2 \leq x' < C_2 a_2$, so $T(a, 0)$ must be the thread $T(C_2-1, 0)$.

*Lemma 9*

The smallest break y in any stride generator satisfies $y = str(T_p) - 1$ or $y = end(T_p)$ for some thread $T_p$ of order p.

(In other words, the smallest break can be found just in front of or at the end of a thread of order p.)

Proof

By condition (C), no thread of order $\leq p+1$ can cross a break y, and so we know that breaks can only arise at the junction of two contiguous threads - say $T_i$ of order i, and $T_j$ of order j. The possibilities are:

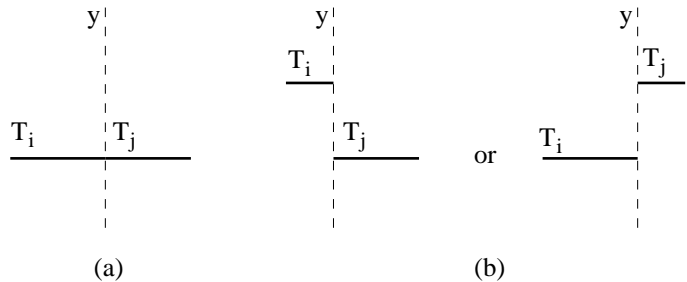

(a)            (b)

Case (a)   (i = j):

In this case, $str(T_i) = str(T_j) - a_2$; so $T_i$ covers $a_2$ values and $T_j$ covers $(a_2-1)$ values. Since $y \geq a_3 - a_2$ (by Lemma 1) this means that all values $0 \leq x < a_3$ are covered by threads of order i. But threads of order 0 are at least as long as threads of order i, and so they, too, must cover the whole stride generator; so $i = p = 0$ by condition (B).

Case (b)   (i != j):

If i = p or j = p the case is proved; so we assume both i, j < p and consider the threads $T_{i+1}$ of order i+1 satisfying $str(T_{i+1}) = str(T_i) - C_1$ and its companion $T_{j+1}$. We know that $len(T_{i+1}) \leq len(T_i)$:

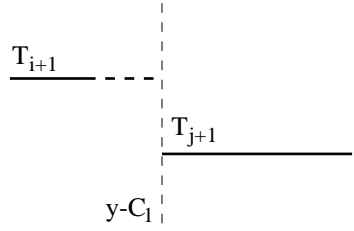

If $T_{i+1}$ does not meet $T_{j+1}$, then there must be a thread V of order k which covers $(y-C_1)$. If k > 0, then V' of order k-1 with $str(V') = str(V) + C_1$ will cross y, so V must be of order zero.

If $T_{i+1}$ meets $T_{j+1}$ and no thread crosses $(y-C_1)$ then $(y-C_1)$ is a break, which contradicts our assumption that y is the smallest break; so in this case, too, some thread V of order zero covers $(y-C_1)$.

Since $y \geq a_3 - a_2$, $(C_2-1)a_2 \leq y-C_1 < C_2 a_2$, and so V must be the thread $T(C_2-1, 0)$:

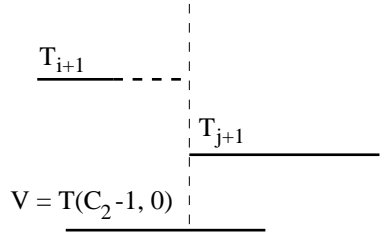

We know that $str(T_{i+1}) < str(T_0)$ because otherwise $T_{i+1}$ would be covered by $T_0$; similarly, $end(T_{j+1}) > end(T_0)$.

Lemma 2 shows that there is a value $x \geq a_3 - a_2$ which is covered only by a thread of order p, and Lemma 8 shows that $(x-C_1)$ is covered only by the thread $T(C_2-1, 0)$; so:

$str(T_{i+1}) < str(T_0) \leq x-C_1 \leq end(T_0) < end(T_{j+1}) \implies str(T_i) < x < end(T_j)$.

But by hypothesis all values in this range are covered by the two threads $T_i$ and $T_j$, both of order < p. So our original assumption that i < p leads to a contradiction, and case (b) is proven.

*Lemma 10*

A stride generator is canonical if one of its breaks is canonical; in other words, either all of the breaks in a stride generator are canonical, or none of them is.

Proof

Using the notation of Lemma 9, suppose y is a break at the junction of two threads $T_i$ and $T_j$ with both i, j < p; then from the proof of that Lemma we know that $y' = y-C_1$ at the junction of the two threads $T_{i+1}, T_{j+1}$ must also be a break. This can only be so if $len(T_{i+1}) = len(T_i)$ which means that $C_2 = 1$. If $C_2 > 1$, there are only two possible positions for breaks in a stride generator: just before, or at the end of, the thread $T_p$.

We first consider $C_2 = 1$, and show that y is a non-canonical break if and only if y' is a non-canonical break:

(i) y is a non-canonical break $\implies$ there is a thread $T_q$ of order q > p+1 which crosses y
$\implies$ there is a thread $T_{q+1}$ which crosses y'
( because $C_2 = 1 \implies len(T_{q+1}) = len(T_q)$ )
$\implies$ y' is a non-canonical break

(ii) y' is a non-canonical break  => there is a thread $T_q$ of order $q > p+1$ which crosses y'
    => there is a thread $T_{q-1}$ which crosses y;
since y is a break, this thread must be of order $> p+1$, and so y is non-canonical.

From this we see that we need only consider the breaks around the thread $T_p$, since other breaks are possible only when $C_2 = 1$ in which case they are all canonical or all non-canonical according to the character of the smallest break(s).

So if there is only one break around $T_p$, the theorem is proved; otherwise the two smallest breaks must be as follows:

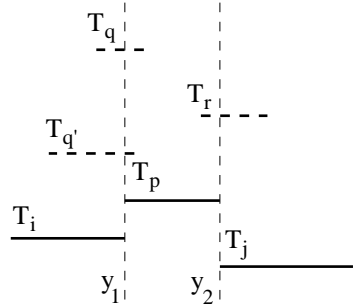

Suppose $y_1$ is non-canonical and so is crossed by $T_q$ for some $q > p+1$; then by similarity let $T_r$ be the thread that is to $T_j$ as $T_q$ is to $T_p$. $T_r$ crosses $y_2$, and so $r > p+1$ - for otherwise $y_2$ would not be a break. So $y_2$ is also non-canonical.

We now apply the argument in the opposite direction to derive the thread $T_{q'}$ that is to $T_i$ as $T_r$ is to $T_p$. Clearly $q' < r < q$, and so by repeated applications we must eventually derive a thread T which crosses $y_1$ and is of order $< p+1$: but this is not possible because $y_1$ is a break.

This contradiction means that our assumption that $y_1$ (or $y_2$) is non-canonical cannot be true: in this configuration, both $y_1$ and $y_2$ are always canonical breaks, and the Lemma is proved.

## 1.3 *The relationship between stride generators and h-bases*

*Lemma 11*

Every h-basis B(A, h) with h-range X has an *underlying* stride generator SG(A, h-k, p) where k is given by $X = (k+1)a_3 + Y$ where $0 \leq Y < a_3 - 1$, and $p \leq k$.

$y = Y+1$ is a break in the stride generator which is either canonical or has break order $> k+1$.

Proof

We first deal with two subsidiary points:

i) We may assume $k \geq 0$ because we are interested only in admissible h-bases.

ii) It is easy to show that Y cannot equal $a_3 - 1$:

Suppose the contrary; this means that $(k+2)a_3$ has no representation, and so $h \leq k+1$. The maximum value that can be represented using at most h values is $ha_3$, and so $X \leq (k+1)a_3$ - which contradicts our assumption that $X = (k+2)a_3 - 1$.

Every value $ka_3 \leq x < (k+1)a_3$ has an h-representation $c_3a_3 + c_2a_2 + c_1$; rewriting, we have:

  $x' + (k-c_3)a_3 = c_2a_2 + c_1$     for  $c_2 + c_1 \leq (h - c_3)$,     $0 \leq x' < a_3$

Writing $i = k - c_3$, we have:

  $x' + ia_3 = c_2a_2 + c_1$             for  $c_2 + c_1 \leq (h - k) + i$,   $0 \leq x' < a_3$,  $k \geq i \geq 0$   (1)

Let p be the smallest value such that (1) is soluble for some $i \leq p$ for all $0 \leq x' < a_3$; then it is clear that conditions (A) and (B) for a stride generator SG(A, h-k, p) are met, with $p \leq k$.

X+1 has no h-representation; writing $y = Y+1$ (and noting that $0 < y < a_3$) we have:

  $y + (k+1)a_3 = c_3a_3 + c_2a_2 + c_1$,  $c_3 + c_2 + c_1 \leq h$, has no solution for $c_1, c_2, c_3 \geq 0$

Writing $j = (k+1) - c_3$, we have:

$y + ja_3 = c_2a_2 + c_1$, $c_2 + c_1 \leq (h - k) + j - 1$, has no solution for $c_1, c_2 \geq 0$, $j \leq k+1$

Since $p \leq k$, this shows condition (C) for a stride generator $SG(A, h-k, p)$ is met; furthermore, y is either a canonical break or has break order $> k+1$.

This correspondence between stride generators and h-bases was first used in [1] where the *potential h-range* $P(h)$ of a stride generator $SG(A, n, p)$ is defined as $P = (h - n + 1)a_3 + y - 1$ where y is its first break (P is called the *potential cover* in [1]). This function is maximised for fixed h to obtain the stride generator $S_{opt} = SG(A_{opt}, n_{opt}, p_{opt})$ with largest potential h-range; it is then shown that $S_{opt}$ is also the stride generator underlying the h-base $B(A_{opt}, h)$, and so P is also the largest h-range that can be realised with any set A. (This is called the *extremal h-range*, and $A_{opt}$ is known as the *extremal basis*; this problem was first solved in 1968 - see [2].)

But the true significance of this Lemma only becomes evident if we suppose that the h-base $B(A, h)$ has the *same* underlying stride generator $SG(A, n, p)$ for all h: that is, n and p are independent of h. If this is so, properties of the h-base which correspond to properties of the stride generator must be independent of h - and it is then straightforward to deduce that $h_1, h_2 \leq h_0$.

It is easy to see that this can only be so if every underlying stride generator is also canonical, a property which was conjectured in [1] but not proved; most of the remainder of this paper is devoted to filling that gap.

## 1.4 Main results

*Theorem 1*

> If $SG(A, n, p)$ has a non-canonical break y with break order q, then $n + q \leq a_2$.
>
> (In fact, it is easy - but tedious - to show $n + q < a_2$; but strict inequality is not necessary for our purposes here.)

Proof

> This is the main new result of this paper, and the proof is given in section 2.

*Theorem 2*

> The stride generator $SG(A, n, p)$ underlying the h-basis $B(A, h)$ is canonical.

Proof

> Suppose this is not the case, and that the stride generator has a break y with break order q.
>
> From Lemma 11, we know that $q > k+1$ and $n = h-k$; so $q > h-n+1 \Rightarrow n+q > h+1$. But for $B(A, h)$ to be admissible we must be able to represent $a_2-1$, and so $h \geq a_2-1$; thus $n+q > a_2$.
>
> This contradicts Theorem 1, and so no such break is possible and the underlying stride generator is canonical as required.

*Theorem 3*

> Let the admissible h-basis $B(A, h)$ have h-range $X(h)$, and $B(A, h+1)$ have h-range $X(h+1)$; then:
>
> $$X(h+1) = X(h) + a_3$$

Proof

> If x has an h-representation, then $(x+a_3)$ has an (h+1)-representation; so we have only to show that there is no (h+1)-representation for $X(h) + a_3 + 1$.
>
> Let $SG(A, n, p)$ be the stride generator underlying $B(A, h)$; we write $X(h) = (k+1)a_3 + Y$, $0 \leq Y < a_3-1$, $y = Y+1$. By Lemma 11 and Theorem 2, $n = h-k$ and y is a canonical break, which means that $y + ja_3 = c_2a_2 + c_1$, $c_2 + c_1 \leq n+j-1$ has no solution for any $j \geq 0$. Writing $j = k + 2 - c_3$, we find $y + (k+2)a_3 = c_3a_3 + c_2a_2 + c_1$, $c_3 + c_2 + c_1 \leq h+1$ has no solution for any $c_3 \leq k+2$: in other words, $X(h) + a_3 + 1$ has no $(h+1)$-representation and the theorem is proved.

Corollory

$h_1 \leq h_0$ for any h-base $A = \{1, a_2, a_3\}$.

The following Lemma - which states that in a canonical stride generator any thread of order $i > p$ is covered by another of order $j \leq p$ - is needed only to prove that $h_2 \leq h_0$.

*Lemma 12*

Let $SG(A, n, p)$ be a canonical stride generator, and let T be a thread of order $i > p$; then there exists a thread U of order $j \leq p$ which covers T.

Proof

By Lemma 5, it is sufficient to show that one thread of order $p+1$ is covered by some thread of order $\leq p$.

Case (a); $p = 0$:

By Lemma 9, the only possible position for a break when $p = 0$ is at the end of a thread of order 0. Since this value cannot be crossed by a thread of order 1, any such thread must be covered by a thread of order 0 and the Lemma is proved.

When $p > 0$, Lemma 9 shows that the stride generator must have a break $y \geq a_3 - a_2$ just in front of or at the end of a thread $T_p$ of order p; we consider these possibilities in turn.

Case (b); $p > 0$, $y = str(T_p) - 1$:

We consider the threads $T_{i+1}$ and $T_{p+1}$ which are displaced $C_1$ to the left of $T_i$ and $T_p$:

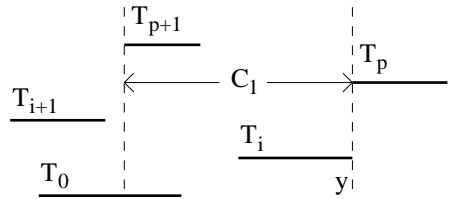

Consider the value $x = str(T_{p+1})$ which must be covered by some thread $T_k$ of order $k \leq p$; we show $k = 0$:

Suppose $k > 0$, and consider the thread $T_{k-1}$ satisfying $str(T_{k-1}) = str(T_k) + C_1$, which covers $str(T_p) = y+1$. $T_{k-1}$ and $T_p$ cannot start in the same position (since then $T_{k-1}$ would cover $T_p$) so $T_{k-1}$ must also cover y, and hence crosses y. This cannot be so since y is a break, so we deduce that $k = 0$.

We now show that if $T_{p+1}$ is not covered by $T_0$, then y is crossed by the thread $T_{i+p+1}$ and so is non-canonical; in other words, if y is canonical then $T_{p+1}$ must be covered by $T_0$ as required.

$str(T_{p+1}) \geq str(T_0)$, and so $T_{p+1}$ is covered by $T_0$ <=> $end(T_{p+1}) \leq end(T_0)$; so we assume $end(T_{p+1}) > end(T_0)$.

$T_p$ is to $T_{i+p+1}$ as $T_0$ is to $T_{i+1}$ so, by similarity:

$str(T_p) - str(T_{i+p+1}) = str(T_0) - str(T_{i+1}) > 0$ (for otherwise $T_{i+1}$ is covered by $T_0$)

=> $str(T_{i+p+1}) \leq y$

$end(T_{i+p+1}) - str(T_p) = end(T_{i+p+1}) - end(T_i) - 1$

$= end(T_{p+1}) - end(T_0) - 1 \geq 0$ (since $T_{p+1}$ is not covered by $T_0$)

=> $end(T_{i+p+1}) \geq y+1$

Case (c); $p > 0$, $y = end(T_p)$:

We consider the threads $T_{j+1}$ and $T_{p+1}$ which are displaced $C_1$ to the left of $T_j$ and $T_p$:

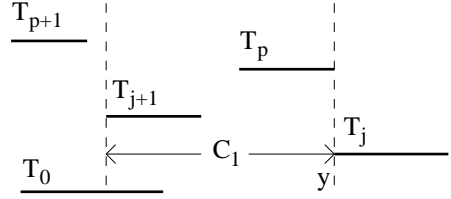

Consider the value $x = str(T_{j+1}) - 1$ which must be covered by some thread $T_k$ of order $k \leq p$; we show $k = 0$:

Suppose $k > 0$, and consider the thread $T_{k-1}$ satisfying $str(T_{k-1}) = str(T_k) + C_1$, which covers $str(T_j) - 1 = y$. $T_{k-1}$ and $T_p$ cannot finish in the same position (since then $T_{k-1}$ would cover $T_p$) so $T_{k-1}$ must also cover $y+1$, and hence crosses $y$. This cannot be so since $y$ is a break, so we deduce that $k = 0$.

We now show that if $T_{p+1}$ is not covered by $T_0$, then $y$ is crossed by the thread $T_{j+p+1}$ and so is non-canonical; in other words, if $y$ is canonical then $T_{p+1}$ must be covered by $T_0$ as required.

$end(T_0) \geq end(T_{p+1})$, and so $T_{p+1}$ is covered by $T_0$ $<=>$ $str(T_{p+1}) \geq str(T_0)$; so we assume $str(T_{p+1}) < str(T_0)$.

$T_p$ is to $T_{j+p+1}$ as $T_0$ is to $T_{j+1}$ so, by similarity:
$end(T_{j+p+1}) - str(T_j) = end(T_{j+p+1}) - end(T_p) - 1$
$= end(T_{j+1}) - end(T_0) - 1 \geq 0$  (for otherwise $T_{j+1}$ is covered by $T_0$)
$=> end(T_{j+p+1}) \geq y+1$
$str(T_j) - str(T_{j+p+1}) = str(T_0) - str(T_{p+1}) > 0$  (since $T_{p+1}$ is not covered by $T_0$)
$=> str(T_{j+p+1}) \leq y$

*Theorem 4*

Let $B(A, h)$ be an admissible h-basis; then:

$x$ has no h-representation $<=>$ $(x+a_3)$ has no (h+1)-representation for all $X(h) < x < ha_3$

Proof

Let $x = (k+r)a_3 + x'$, $0 \leq x' < a_3$, $r \geq 1$; $n = h - k$ as usual.

Then $x$ has no h-representation means $x' + (k+r)a_3 = c_3 a_3 + c_2 a_2 + c_1$, $c_3 + c_2 + c_1 \leq h$ has no solution for $0 \leq c_3 \leq k+r$.

Writing $i = (k+r-c_3)$ we have:
$$x' + ia_3 = c_2 a_2 + c_1, \quad c_2 + c_1 \leq n + i - r \quad (1)$$
has no solution for $0 \leq i \leq k+r$. Similarly, $(x+a_3)$ has no (h+1)-representation means that (1) has no solution for $0 \leq i \leq k+r+1$.

Solutions to (1) can be found by taking the thread diagram for the underlying stride generator $SG(A, n, p)$ and reducing the length of each thread by $r$; (1) has a solution if and only if there is a truncated thread of order $i$ which covers $x'$. $SG(A, n, p)$ is canonical by Theorem 2, and so by Lemma 12 we need only consider threads of order $\leq p$; this means that if (1) has no solution for $i \leq p$, then it has no solution at all.

Since $p \leq k$ by Lemma 11, the theorem is proved.

Corollory

$h_2 \leq h_0$ for any h-base $A = \{1, a_2, a_3\}$.

## 2  Every non-canonical stride generator has $n + q \leq a_2$

### 2.1  *Preparatory remarks*

Before outlining the proof of Theorem 1, we require a few more definitions and lemmas.

*Lemma 13*

> If SG(A, n, p) is a canonical stride generator, then no stride generator SG(A, n', p') exists for $n' < n$.

Proof

> The thread diagram for SG(A, n', p') is obtained from that of SG(A, n, p) by reducing the length of each thread by (n - n'). This 'uncovers' any canonical break y in SG(A, n, p), thus showing that y has no n'-generation for any $n' < n$.

We say that SG(A, n, p) is the *fundamental* stride generator for A if there is no other stride generator SG(A, n', p') with $n' > n$.

> It is easy to see that the fundamental stride generator SG(A, $n_1$, $p_1$) is the first in a series of stride generators SG(A, $n_i$, $p_i$) with $n_{i+1} < n_i$, $p_{i+1} > p_i$ that terminates with a canonical stride generator SG(A, $n_t$, $p_t$); each stride generator for $i < t$ is non-canonical. If the fundamental stride generator is canonical, then t = 1. These stride generators are the only stride generators SG(A, n, p) for the set A.

Our proof of Theorem 1 proceeds as follows:

> We first show that any non-canonical fundamental stride generator has order $p \geq 2$, and that its thread diagram in the range $(C_2 - 1)a_2 \leq x < C_2 a_2$ has a particular form: it has the appearance of either an ascending ($C_1 > a_2/2$) or descending ($C_1 < a_2/2$) staircase.

> For $C_2 > 1$, we determine an upper bound $q_{max}$ such that no thread T of order $q_{max}$ exists within this range for the fundamental stride generator. This means that $q < q_{max}$ for any break y in the fundamental stride generator (or in any derived from it). We then show that $n + (q_{max} - 1) \leq a_2$, which proves the result.

> A different approach is necessary when $C_2 = 1$. In this case we determine the upper bound $q_{max}$ by demonstrating the existence of a thread of order $q_{max}$ that is covered by $T_0 = T(0, 0)$; we know by Lemma 5 that this means that all threads of order $\geq q_{max}$ are covered by threads of order $< q_{max}$, and so $q < q_{max}$ as before.

### 2.2  *The form of fundamental stride generators*

*Lemma 14*

> The fundamental stride generator SG(A, n, p) for a set A is canonical if it is of order 0 or 1.
>
> Otherwise $p \geq 2$ and SG(A, n, p) has a thread diagram whose format in the range $(C_2 - 1)a_2 \leq x < C_2 a_2$ corresponds to one of the four possibilities shown below:
>
>> In cases (A2) and (D2), the stride generator is canonical.
>>
>> In cases (A1) and (D1) it may or may not be canonical.
>
> Note: (A1) and (A2) are characterised by $str(T_{i+1}) > str(T_i)$ for $0 < i < p$, and
>
>> for (A1):  $y = end(T_{p-1}) = str(T_p) - 1$;     for (A2):  $y = end(T_p) = str(T'_0) - 1$.
>
> (D1) and (D2) are characterised by $str(T_{i+1}) < str(T_i)$ for $0 < i < p$, and
>
>> for (D1):  $y = end(T_p) = str(T_{p-1}) - 1$;     for (D2):  $y = end(T_0) = str(T_p) - 1$.

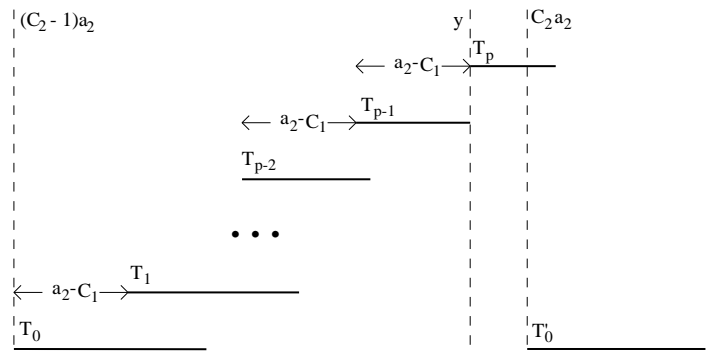
A1: Ascending staircase, type 1

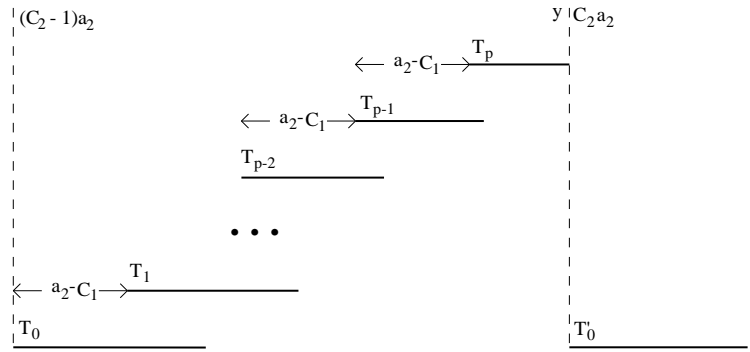
A2: Ascending staircase, type 2

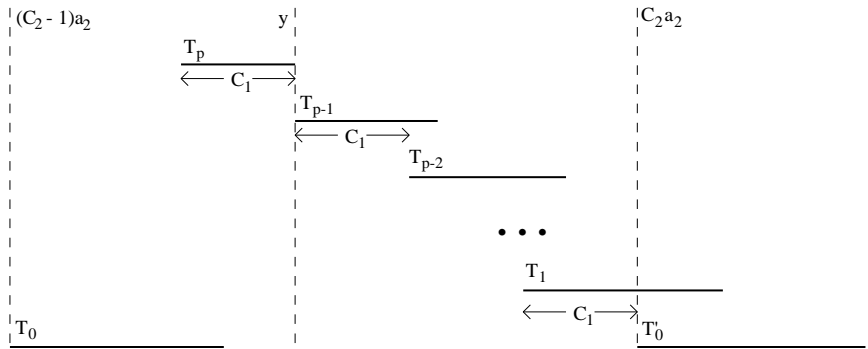
D1: Descending staircase, type 1

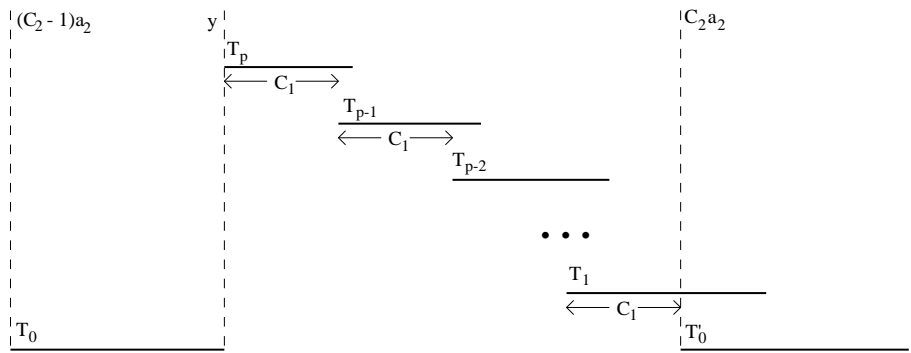
D2: Descending staircase, type 2

*Page 14*

Proof

We first note some properties of these thread diagrams.

If every value $(C_2 - 1)a_2 \leq x < C_2 a_2$ is covered by some thread, then so are all values $0 \leq x < (C_2 - 1)a_2$. Furthermore, all values $C_2 a_2 \leq x < a_3$ are also covered provided that any break $y \geq a_3 - a_2$. So to show that such a thread diagram corresponds to a stride generator we need only consider threads in the range $(C_2 - 1)a_2 \leq x < C_2 a_2$ provided that we show also that the smallest break $y \geq a_3 - a_2$.

In descending staircases (D1 and D2 above):

$\text{str}(T_{i+1}) = \text{str}(T_i) - C_1$ for $i \geq 1$; $\qquad \text{str}(T_1) = C_2 a_2 - C_1$

$\text{len}(T_{i+1}) = \text{len}(T_i) - (C_2 - 1)$ for $i \geq 1$; $\qquad \text{len}(T_1) = \text{len}(T_0) - C_2$

$\Rightarrow \text{len}(T_0) > \text{len}(T_1) \geq \text{len}(T_2) \geq \ldots \geq \text{len}(T_{p-1}) \geq \text{len}(T_p)$

In ascending staircases (A1 and A2 above):

$\text{str}(T_{i+1}) = \text{str}(T_i) + (a_2 - C_1)$ for $i \geq 0$

$\text{len}(T_{i+1}) = \text{len}(T_i) - C_2$ for $i \geq 0$

$\Rightarrow \text{len}(T_0) > \text{len}(T_1) > \text{len}(T_2) > \ldots > \text{len}(T_{p-1}) > \text{len}(T_p)$

The largest value of n that makes sense to consider is that which causes $T_0 = T(C_2 - 1, 0)$ to cover this entire range. The only possible position for a break y is at the end of $T_0$, and $y = C_2 a_2 - 1 \geq a_3 - a_2$ as required. If $T_1$ does not cross y, this is a zero order stride generator and hence the fundamental stride generator for A:

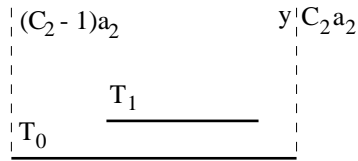

We see that in this case $T_1$ is covered by $T_0$, and so the stride generator is canonical by Lemma 5.

Now suppose that $T_1$ crosses the end of $T_0$; we reduce n until the two threads together *just* cover the range. We consider three separate cases: $C_1 = a_2/2$, $C_1 < a_2/2$ and $C_1 > a_2/2$.

*When $C_1 = a_2/2$* the only possible thread arrangement is:

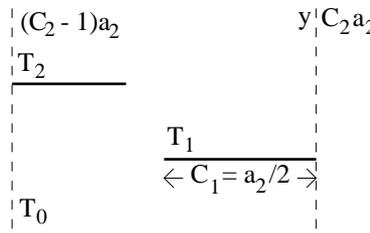

since $\text{len}(T_2) \leq \text{len}(T_1) < \text{len}(T_0)$. This has a break $y = C_2 a_2 - 1$, and so $y \geq a_3 - a_2$ as required. Since $T_2$ is covered by $T_0$, this fundamental stride generator of order 1 is canonical by Lemma 5.

*When $C_1 < a_2/2$*, there are two possibilities:

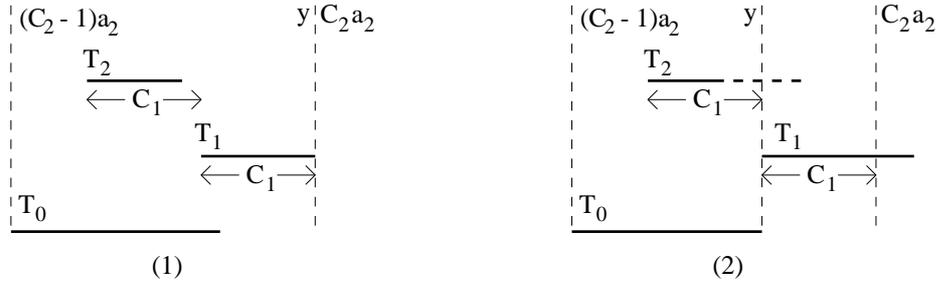

(1)                          (2)

In case (1), we have a canonical fundamental stride generator of order 1, because:

$y = C_2 a_2 - 1 \leq a_3 - a_2$

$T_2$ is covered by $T_0$ because:

$str(T_2) = C_2 a_2 - 2C_1 > (C_2 - 1)a_2 = str(T_0)$     since $C_1 < a_2/2$

$end(T_2) = end(T_1) - C_1 - (C_2 - 1) \leq str(T_1) - 1 \leq end(T_0)$

In case (2), we have a canonical fundamental stride generator of order 1 if $T_2$ is covered by $T_0$, since $y = C_2 a_2 - C_1 - 1 \geq a_3 - a_2$ when $C_1 < a_2/2$; otherwise $T_2$ crosses $y$ and we have the beginning of a descending staircase.

Once again, we reduce n until the threads $T_0$, $T_2$ and $T_1$ just cover the range: the result is one of the two possibilities illustrated below with $k = 2$ (note that it is impossible for $end(T_1) = C_2 a_2$ while $end(T_2) > str(T_1)$ because $len(T_1) \geq len(T_2)$):

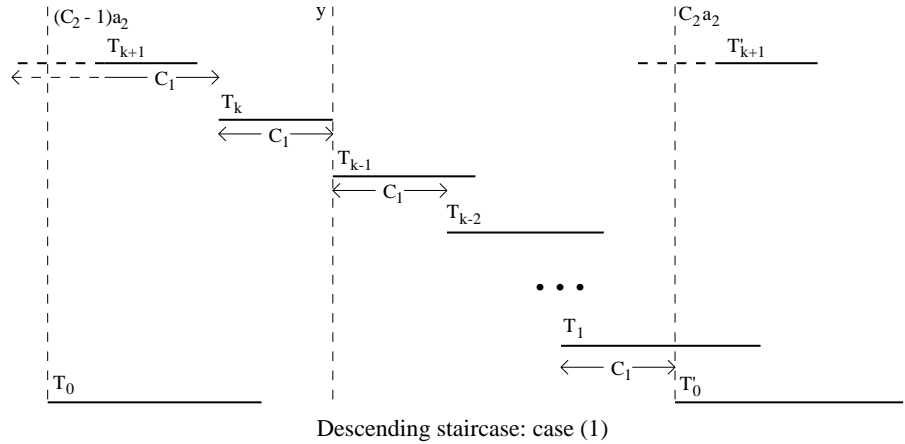

Descending staircase: case (1)

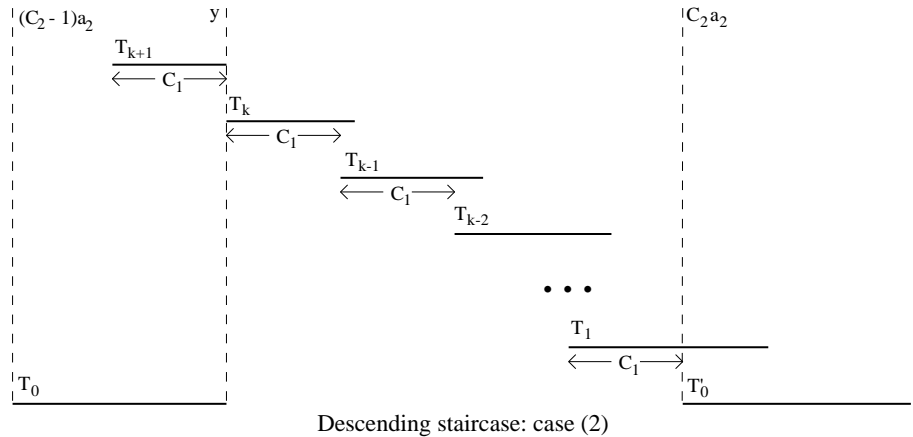

Descending staircase: case (2)

We now show in general for $k \geq 2$:

(A) Case (1) describes a fundamental stride generator of order k which may or may not be canonical.

(B) If $T_{k+1}$ does not exist or is covered by $T_0$, case (2) describes a canonical fundamental stride generator of order k.

(C) If $T_{k+1}$ crosses the end of $T_0$, case (2) does not describe a stride generator and we reduce n further until case (1) or case (2) for $k' = k + 1$ arises.

So to find the fundamental stride generator for some $C_1 < a_2/2$ we repeat this process for k = 2, 3, ... until one of cases (A) or (B) arises - and the fundamental stride generator can only be non-canonical in case (A). Note that this procedure must terminate because case (1), or case (2) where $T_{k+1}$ does not exist, will eventually arise.

We first note that in both cases:

$$y \geq \text{end}(T_0) = \text{str}(T_0) + \text{len}(T_0) - 1 \geq \text{str}(T_0) + \text{len}(T_k) \quad \text{since } \text{len}(T_0) > \text{len}(T_k)$$
$$\geq \text{str}(T_0) + C_1 = a_3 - a_2.$$

Next we note that if $T_{k+1}$ does not exist, then both cases (1) and (2) describe a canonical stride generator of order k; so now we assume that $\text{len}(T_{k+1}) > 0$.

In case (1), we know that no thread of order k+1 can cross y:

$T_{k+1}$ cannot cross y because $\text{end}(T_{k+1}) < \text{str}(T_k) \leq y$.

$T'_{k+1}$ cannot cross y because:

$$\text{str}(T'_{k+1}) = \text{str}(T_{k+1}) + a_2 = \text{str}(T_k) - C_1 + a_2 > \text{str}(T_k) - C_1 + 2C_1$$
$$> \text{str}(T_k) + C_1 - 1 = y$$

So case (1) represents a stride generator of order k which may or may not be canonical.

In case (2) we know that $\text{str}(T_{k+1}) > \text{str}(T_0)$ because:

$T_k$ at least meets $T_{k-1}$ => $\text{len}(T_k) \geq C_1$ => $\text{len}(T_0) > C_1$, and so
$\text{str}(T_{k+1}) = \text{str}(T_k) - C_1 = \text{end}(T_0) + 1 - C_1 > \text{str}(T_0)$.

So if $\text{end}(T_{k+1}) \leq \text{end}(T_0) = y$, $T_{k+1}$ is covered by $T_0$ and case (2) describes a canonical stride generator of order k.

If $\text{end}(T_{k+1}) > \text{end}(T_0) = y$, $T_{k+1}$ crosses y and so case (2) does not describe a stride generator at all, and we must reduce n to 'reveal' $T_{k+1}$ until the threads $T_0, T_{k+1}, T_k, ... T_1$ just cover the range. Since $\text{len}(T_{k+1}) \leq \text{len}(T_k) \leq ... \leq \text{len}(T_1) < \text{len}(T_0)$ we know that this procedure will result in case (1) or case (2) where k is replaced by k+1 throughout.

*When $C_1 > a_2/2$*, there is only one possible arrangement - see (1) below; this is because two contiguous threads $T_0$ and $T_1$ cannot cover the range - as shown in (2):

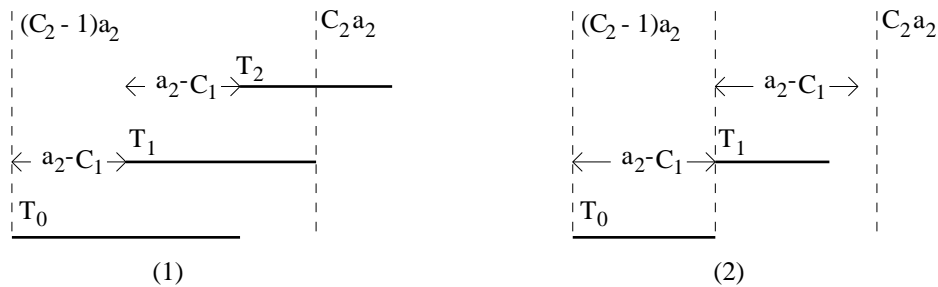

(1)                         (2)

By similarity, $T_2$ must cross $y = C_2 a_2 - 1$, and so this is not a stride generator; instead it is the beginning of an ascending staircase. So we reduce n until the threads $T_0$, $T_1$ and $T_2$ just cover the range, resulting in one of the possibilities shown below with k = 2:

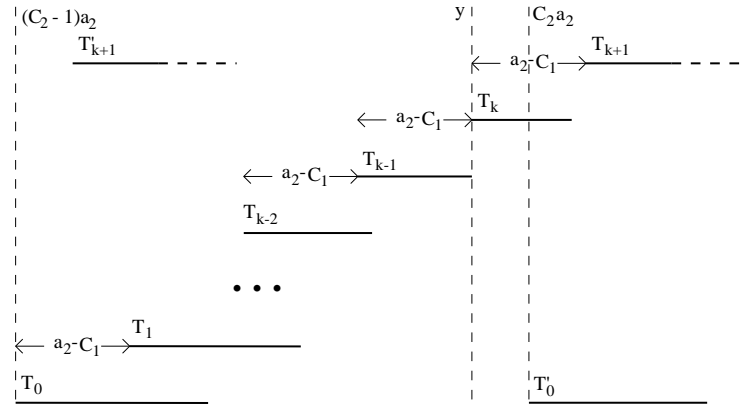

Ascending staircase: case (1)

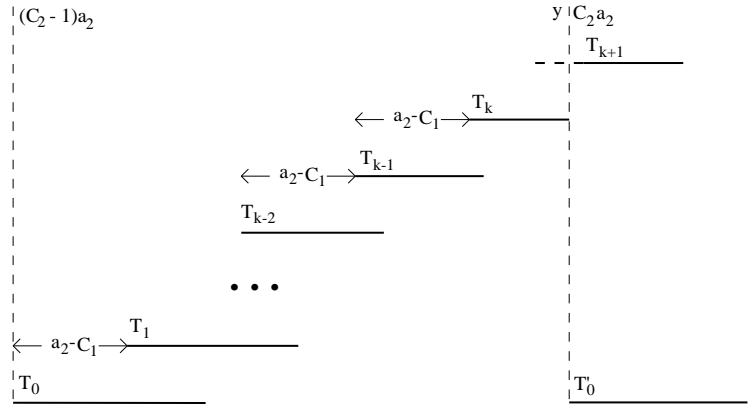

Ascending staircase: case (2)

We now show in general for $k \geq 2$:

(A) Case (1) describes a fundamental stride generator of order k which may or may not be canonical.

(B) If $T_{k+1}$ does not exist, or is covered by $T'_0$, case (2) describes a canonical fundamental stride generator of order k.

(C) If $T_{k+1}$ crosses the end of $T_k$, case (2) does not describe a stride generator and we reduce n further until case (1) or case (2) for $k' = k+1$ arises.

So to find the fundamental stride generator for some $C_1 > a_2/2$ we repeat this process for $k = 2, 3, \ldots$ until one of cases (A) or (B) arises - and the fundamental stride generator can only be non-canonical in case (A). This procedure must terminate because case (1), or case (2) where $T_{k+1}$ does not exist, will eventually arise.

We first note that if $T_{k+1}$ does not exist, then both cases (1) and (2) describe a canonical stride generator of order k; so we now assume $len(T_{k+1}) > 0$.

In case (1) we know that no thread of order k+1 can cross y:

$T_{k+1}$ cannot cross y because $str(T_{k+1}) > end(T_k) > y$.

$T'_{k+1}$ cannot cross y because:

$$\begin{aligned} end(T'_{k+1}) &= end(T_{k+1}) - a_2 + 1 < end(T_k) + a_2 - C_1 - a_2 + 1 \\ &\leq str(T_k) + len(T_k) - C_1 < str(T_k) + a_2 - 2C_1 \\ &= y + 1 + a_2 - 2C_1 \leq y \quad (\text{since } a_2 < 2C_1) \end{aligned}$$

Furthermore, $y \geq C_2 a_2 - (len(T_k) + 1) \geq C_2 a_2 - (a_2 - C_1) = a_3 - a_2$.

So case (1) represents a stride generator of order k which may or may not be canonical.

In case (2) we know that $end(T_{k+1}) < end(T'_0)$ because:

$$end(T_{k+1}) < end(T_k) + (a_2 - C_1) = (C_2 a_2 - 1) + (a_2 - C_1), \text{ and}$$

$$\text{end}(T'_0) = \text{end}(T_0) + a_2 - 1 \geq (C_2 - 1)a_2 + a_2 - C_1 + a_2 - 1 = (C_2 a_2 - 1) + (a_2 - C_1)$$

Furthermore, $y = C_2 a_2 - 1 \geq a_3 - a_2$.

So if $\text{str}(T_{k+1}) \geq C_2 a_2$ then $T_{k+1}$ is covered by $T'_0$ and case (2) decribes a canonical stride generator of order k.

If $\text{str}(T_{k+1}) \leq C_2 a_2 - 1 = \text{end}(T_k) = y$, $T_{k+1}$ crosses y and so case (2) does not describe a stride generator at all, and we must reduce n to 'reveal' $T_{k+1}$ until the threads $T_0, T_1, \ldots T_{k+1}$ together just cover the range. Since $\text{len}(T_{k+1}) < \text{len}(T_k) < \ldots < \text{len}(T_1) < \text{len}(T_0)$ we know that this process will result in case (1) or case (2) where k is replaced by k+1 throughout.

This completes the proof of Lemma 14.

We may consider how the order p of the fundamental stride generator varies as $C_1$ varies for a fixed (but small) value of $C_2$. The proof of Lemma 14 shows that:

For $C_1 > a_2/2$, p increases from 2 as $C_1$ increases towards some critical value $X_0$; for all values $X_0 \leq C_1 < a_2$ (and for $C_1 = 0$), p = 0.

For $1 \leq C_1 \leq a_2/2$, p decreases from some large value as $C_1$ increases towards some critical value $X_1$; for all $X_1 \leq C_1 \leq a_2/2$, p = 1.

*Lemma 15*

If $a_2 - C_2 \leq C_1 < a_2$, or if $C_1 = 0$, the fundamental stride generator for A is of order 0.

Proof

The critical part of the thread diagram for a zero order stride generator when $C_1 > 0$ has the following appearance:

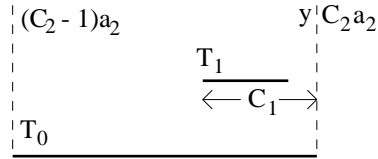

We require $\text{len}(T_1) \leq C_1$; but $\text{len}(T_1) = \text{len}(T_0) - C_2 = a_2 - C_2$; so $C_1 \geq a_2 - C_2$ as required.

In the special case of $C_1 = 0$, $\text{str}(T_1) = \text{str}(T_0)$ and $\text{len}(T_1) \leq \text{len}(T_0)$; so $T_1$ is always covered by $T_0$ and the stride generator is of order 0.

*Lemma 16*

If $a_2 \geq 2C_1 \geq a_2 - 2C_2 + 1$, then the fundamental stride generator for A is of order 0 or 1.

Proof

The proof of Lemma 14 shows that if the fundamental stride generator is not of order 0, then it is of order 1 if $C_1 = a_2/2$, and may be of order 1 for $C_1 < a_2/2$ in the following situations:

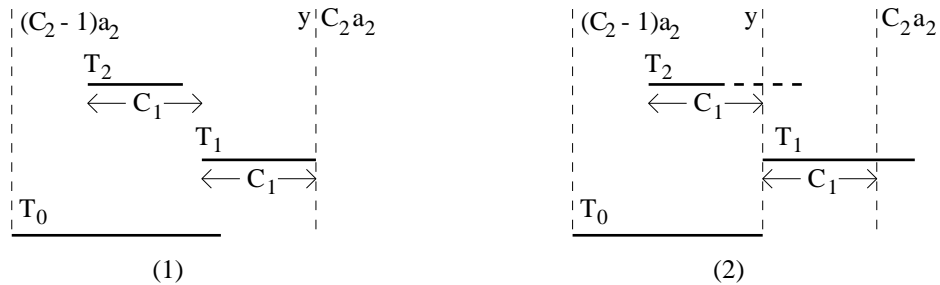

(1)                              (2)

In both cases, $\text{len}(T_0) = n - C_2 + 2$, $\text{len}(T_1) = n - 2C_2 + 2$ and $\text{len}(T_2) = n - 3C_2 + 3$.

In case (1):

$\text{len}(T_1) = C_1 \Rightarrow n = C_1 + 2C_2 - 2$

$T_0$ must at least meet $T_1$, so:

$\text{len}(T_0) \geq a_2 - C_1 \Rightarrow C_1 + C_2 \geq a_2 - C_1 \Rightarrow 2C_1 \geq a_2 - C_2 \geq a_2 - 2C_2 + 1$

In case (2):

$\text{len}(T_0) = a_2 - C_1 \Rightarrow n = a_2 - C_1 + C_2 - 2$

$T_2$ must not cross the end of $T_0$, so:

$\text{len}(T_2) \leq C_1 \Rightarrow a_2 - C_1 - 2C_2 + 1 \leq C_1 \Rightarrow 2C_1 \geq a_2 - 2C_2 + 1$

## 2.3 The descending staircase:  $C_1 < a_2/2$

We know from Lemma 14 that any non-canonical fundamental stride generator with $C_1 < a_2/2$ has the following form of thread diagram for some $p \geq 2$:

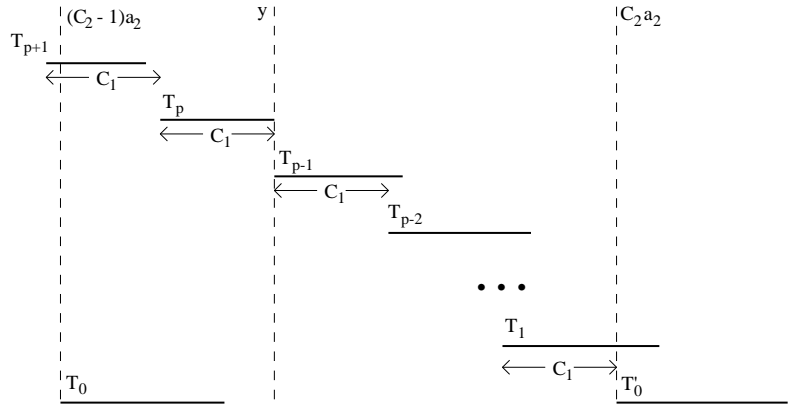

Descending staircase: case (1)

We know that $T_{p+1}$ crosses $\text{str}(T_0) - 1$ because:

$\text{str}(T_{p+1}) < \text{str}(T_0)$ because:

$\text{end}(T_{p+1}) \leq \text{end}(T_0)$, so $T_{p+1}$ is covered by $T_0$ if $\text{str}(T_{p+1}) \geq \text{str}(T_0)$ and the stride generator would then be canonical.

$\text{end}(T_{p+1}) \geq \text{str}(T_0)$ because:

$\text{end}(T_{p+1}) = \text{end}(T_p) - C_1 - (C_2 - 1) > \text{end}(T_0) - C_1 - (C_2 - 1)$, so
$\text{end}(T_{p+1}) - \text{str}(T_0) > \text{len}(T_0) - C_1 - C_2$; but $\text{len}(T_1) = \text{len}(T_0) - C_2 \geq C_1$, so
$\text{end}(T_{p+1}) - \text{str}(T_0) > 0$.

### 2.3.1 General bounds

We start with an improved bound on $C_1$ obtained immediately from Lemma 16:

$$2C_1 \leq a_2 - 2C_2 \qquad - (0)$$

We obtain bounds for $C_1$ and a formula for n as follows:

We see immediately that $(p+1)C_1 > a_2 > pC_1$. Here we need only the lower bound for $C_1$, because we develop a better upper bound below in (3):

$$(p+1)C_1 > a_2 \qquad - (1)$$

Next, we note:

$T'_0 = T(C_2, 0) \Rightarrow T_1 = T(2C_2, 1) \Rightarrow ... \Rightarrow T_p = T((p+1)C_2, p)$

so   $\text{len}(T_p) = C_1 = n + p - (p+1)C_2 + 1$

$\Rightarrow \qquad\qquad\qquad\qquad\qquad n = C_1 + (p + 1)(C_2 - 1) \qquad\qquad\qquad - (2)$

We have:

$\text{end}(T_0) = \text{str}(T_0) + n - (C_2 - 1) = \text{str}(T_0) + C_1 + p(C_2 - 1)$   by (2)

$\text{end}(T_p) = \text{str}(T_0) + a_2 - (p - 1)C_1 - 1$

Now $\text{end}(T_0) < \text{end}(T_p)$, so $C_1 + p(C_2 - 1) < a_2 - (p - 1)C_1 - 1$

$\Rightarrow \qquad\qquad\qquad\qquad\qquad p(C_1 + C_2 - 1) < a_2 - 1 \qquad\qquad\qquad - (3)$

We can now derive an upper bound for $C_2$ as follows:

From (3), $pC_2 < a_2 - pC_1 + p - 1$

From (1), $pC_1 > pa_2/(p+1)$, so $pC_2 < a_2/(p+1) + p - 1$

$\Rightarrow \qquad\qquad\qquad\qquad\qquad C_2 < a_2/(p(p+1)) + 1 \qquad\qquad\qquad - (4)$

### 2.3.2  *The case for $C_2 \geq 2$*

In this section, we assume $C_2 \geq 2$.

We obtain an upper bound for q as follows:

Let $T_j$ be the thread of order j such that $(C_2 - 1)a_2 \leq \text{str}(T_j) < C_2 a_2$:

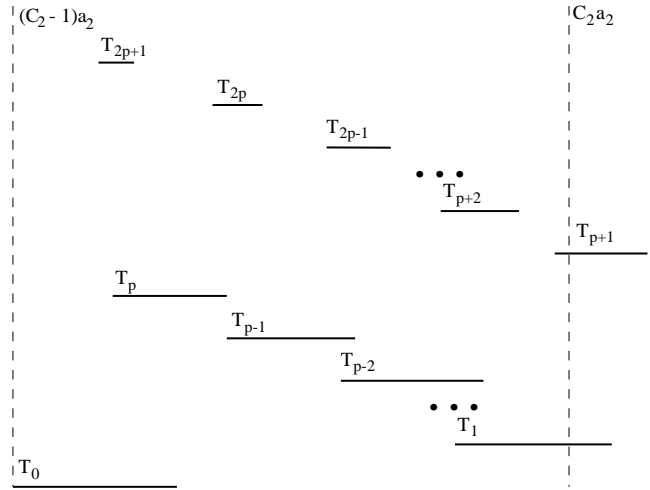

If $T_{j+1}$ and $T_j$ are two steps on the same staircase, $\text{len}(T_{j+1}) = \text{len}(T_j) - (C_2 - 1)$; if they are at opposite ends of the range (eg j = p), $\text{len}(T_{j+1}) = \text{len}(T_j) - C_2$; so $\text{len}(T_{j+1}) \leq \text{len}(T_j) - (C_2 - 1)$ for $j \geq 0$. Since $\text{len}(T_1) = \text{len}(T_0) - C_2 = n - 2C_2 + 2$, we have:

$\text{len}(T_j) \leq n - (j + 1)(C_2 - 1)$

Suppose Q is the smallest value of j such that $\text{len}(T_j) \leq 0$; that is, $T_{Q-1}$ is the highest order thread that is present in the range. If the fundamental stride generator $S = SG(A, n, p)$ represented by this thread diagram is non-canonical, then the break order q of any break must be less than Q. Furthermore, the same must be true for any non-canonical stride generator $S' = SG(A, n', p')$ derived from S, since the thread diagram for S' is derived from that for S by removing (n - n') units from the end of each thread.

Now $\text{len}(T_Q) \leq 0$ if $n - (Q + 1)(C_2 - 1) \leq 0 \Rightarrow Q \geq n/(C_2 - 1) - 1$; so:

$\qquad\qquad\qquad\qquad\qquad q < n/(C_2 - 1) - 1 \qquad\qquad\qquad - (5)$

We can now derive an upper bound for n + q:

From (5): $n + q < n(1 + 1/(C_2 - 1)) - 1 = n(C_2/(C_2 - 1)) - 1$

From (3): $pC_1 < a_2 - 1 - p(C_2 - 1) \Rightarrow C_1 < (a_2 - 1)/p - (C_2 - 1)$

Substituting for $C_1$ in (2) gives: $n < (a_2 - 1)/p + p(C_2 - 1)$

So:
$$n + q < (\,(a_2 - 1)/p + p(C_2 - 1)\,)(C_2/(C_2 - 1)) - 1 \quad - (6)$$

Hence $n + q < (a_2/p + p(C_2 - 1))\,(C_2/(C_2 - 1))$, and substituting for the first occurrence of $(C_2 - 1)$ using (4) gives:
$$n + q < (a_2/p + a_2/(p+1))\,(C_2/(C_2 - 1)) \quad - (7)$$

We are now ready to prove that $n + q \leq a_2$; we take the cases $p \geq 4$, $p = 3$ and $p = 2$ separately.

### When $p \geq 4$:

From (7) we have:

$\quad n + q < 2(\,(2p+1)/(p(p+1))\,)a_2 \qquad$ since $C_2 \geq 2$

$\Rightarrow n + q < (9/10)a_2 < a_2 \qquad$ for all $p \geq 4$

### When $p = 3$:

From (7) we have:

$\quad n + q < (a_2/3 + a_2/4)(C_2/(C_2 - 1)) = (7/12)a_2\,(C_2/(C_2 -1)) < a_2$ for $C_2 \geq 3$

This leaves $C_2 = 2$; we substitute directly in (6):

$\quad n + q < 2((a_2 - 1)/3 + 3) - 1 = (2a_2 + 13)/3$, which is $\leq a_2$ so long as $a_2 \geq 13$.

From (4) we find $2 < a_2/12 + 1 \Rightarrow a_2 > 12$, which is just sufficient.

### When $p = 2$:

From (7) we have:

$\quad n + q < (a_2/2 + a_2/3)(C_2/(C_2 - 1)) = (5/6)a_2\,(C_2/(C_2 - 1)) \leq a_2$ for $C_2 \geq 6$

This leaves $C_2 = 2, 3, 4$ and $5$ to be considered.

For $C_2 = 5$, we substitute directly in (6):

$\quad n + q < ((a_2 - 1)/2 + 8)(5/4) - 1 = (5a_2 + 67)/8$, which is $\leq a_2$ so long as $a_2 \geq 23$.

From (4) we find $5 < a_2/6 + 1 \Rightarrow a_2 > 24$, which is sufficient.

For $C_2 = 4$, we similarly substitute directly in (6):

$\quad n + q < ((a_2 - 1)/2 + 6)(4/3) - 1 = (4a_2 + 38)/6$, which is $\leq a_2$ so long as $a_2 \geq 19$.

From (4) we find $4 < a_2/6 + 1 \Rightarrow a_2 > 18$, which is just sufficient.

For $C_2 = 3$ a different approach is necessary:

$\quad$ (2) gives: $n = C_1 + 6$

$\quad$ (5) gives: $q < n/2 - 1$, so $n + q < (3/2)C_1 + 8$

$\quad$ (0) gives: $2C_1 \leq a_2 - 6$, so $n + q < (3a_2 + 14)/4$, which is $\leq a_2$ so long as $a_2 \geq 14$.

From (4) we find $3 < a_2/6 + 1 \Rightarrow a_2 > 12$; this leaves $a_2 = 13$ to consider in more detail:

$\quad$ (0) gives $2C_1 \leq 7 \Rightarrow C_1 \leq 3$, and (1) gives $3C_1 > 13 \Rightarrow C_1 \geq 5$; so there is no such case to consider after all.

For $C_2 = 2$ we proceed in a similar way:

    (2) gives: $n = C_1 + 3$

    (5) gives: $q < n - 1$, so $n + q < 2C_1 + 5 \Rightarrow n + q \leq 2C_1 + 4$

    (0) gives: $2C_1 \leq a_2 - 4$, so $n + q \leq a_2$ as required.

(This is one of the two cases where we prove only that $n + q \leq a_2$. To obtain strict inequality we have to use an improved upper bound for q which takes account of the extra reduction in the length of the threads $T_j$ that happens each time a new staircase starts (cf (5) in section 2.4.2 below). Even then, there remain six explicit stride generators which have to be shown individually to satisfy $n + q < a_2$.)

(The details are as follows. We know that:
    $\text{len}(T_0) = n + 1$
    $\text{len}(T_1) = \text{len}(T_0) - 2$
    $\text{len}(T_2) = \text{len}(T_1) - 1$
    $\text{len}(T_3) = \text{len}(T_2) - 2$
    $\text{len}(T_4) = \text{len}(T_3) - 1$
    $\text{len}(T_5) = \text{len}(T_4) - 1$ or $\text{len}(T_4) - 2$
    etc.

We deduce that $\text{len}(T_j) \leq n + 1 - (4/3)j$, so that $\text{len}(T_j) \leq 0$ as soon as $j \geq 3(n+1)/4$. So $q < 3(n+1)/4$, and we have:
    $n + q < (7n+3)/4$

Substituting $n = C_1 + 3$, we have $n + q < (7C_1 + 24)/4$, and substituting $2C_1 \leq a_2 - 4$ gives $n + q < (7a_2 + 20)/8$; so $n + q < a_2$ when $a_2 \geq 20$. From (4) - or just before - we have $C_2 < a_2/6 + 1/2$, which gives $a_2 > 9$; so we have only to consider $a_2$ from 10 to 19 inclusive. The following table gives for each $a_2$ in the range:
    X - the smallest value of $C_1$ which makes $(7C_1 + 24)/4 > a_2$
    Y - the largest value of $C_1$ which allows $p \geq 2$, and hence:
    Z - the values of $C_1$ which we must consider
    n
    $q_{max}$ - the largest value of q satisfying both $q < n - 1$ and $q < 3(n+1)/4$ ($q < n + 1$ is a better bound only for $a_2 = 10$)

| $a_2$ | 10 | 11 | 12 | 13 | 14 | 15 | 16 | 17 | 18 | 19 |
|---|---|---|---|---|---|---|---|---|---|---|
| X | 3 | 3 | 4 | 5 | 5 | 6 | 6 | 7 | 7 | 8 |
| Y | 3 | 3 | 4 | 4 | 5 | 5 | 6 | 6 | 7 | 7 |
| Z | 3 | 3 | 4 | - | 5 | - | 6 | - | 7 | - |
| n | 6 | 6 | 7 |   | 8 |   | 9 |   | 10 |   |
| $q_{max}$ | 4 | 5 | 5 |   | 6 |   | 7 |   | 8 |   |

We see that in all the cases that we must consider, $n + q_{max} = a_2$. Examination of the individual stride generators shows that the first three ( {1, 10, 23}, {1, 11, 25}, {1, 12, 28} ) are all canonical, and the last three ( {1, 14, 33}, { 1, 16, 38}, {1, 18, 43} ) are non-canonical with $q = 4$; so we see that $n + q < a_2$ in all cases.)

This completes the proof for the case $C_1 < a_2/2$, $C_2 \geq 2$.

### 2.3.3 *The case for $C_2 = 1$*

From (2) we have $n = C_1$, and so:

$$n = C_1 < a_2/2 \quad\quad\quad - (8)$$

and, since $pC_1 < a_2 < (p+1)C_1$, we can write:

$$a_2 = pn + s \text{ where } 1 \leq s < n \quad\quad\quad - (9)$$

As before, let $T_j$ be the thread of order j such that $(C_2 - 1)a_2 = 0 \leq \text{str}(T_j) < C_2 a_2 = a_2$. We now consider those threads $S_i$ which satisfy $0 \leq \text{str}(S_i) < n$:

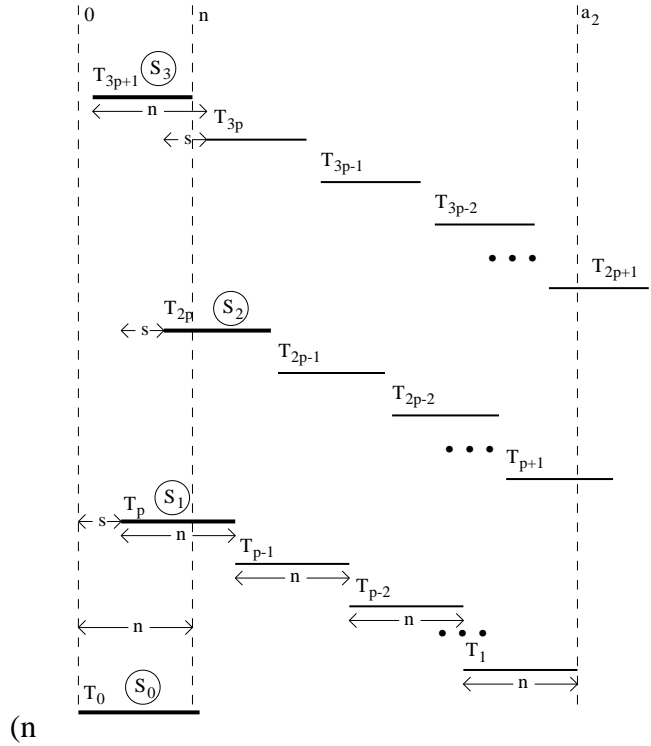

Since $\text{len}(T_0) = n+1$, each thread $S_i$ at least satisfies $\text{str}(T_0) \le \text{str}(S_i) < \text{end}(T_0)$. We will show below that there is always a thread $S_i$ that satisfies $\text{end}(S_i) \le \text{end}(T_0)$ and so is covered by $T_0$, thus providing an upper bound for $q$; it turns out that this bound is sufficient to show that $n + q < a_2$.

First we derive formulae for the order, position and length of thread $S_i$, by observing that $S_i$ is derived from the thread $X = T(i(p+1), ip)$ as follows:

We have $\text{str}(X) = i(p+1)a_2 - ipa_3 = i(a_2 - pn) = is$ and $\text{len}(X) = n + ip - i(p+1) + 1 = (n+1) - i$.

With $C_2 = 1$, a thread $X_1 = T(c, i)$ implies the existence of a further thread $X_2 = T(c+1, i+1)$ of the same length where $\text{str}(X_2) = \text{str}(X_1) - C_1 = \text{str}(X_1) - n$. Let us write:

$$is = kn + t \text{ where } 0 \le t < n \qquad - (10)$$

Then from $X$ we derive thread $Y = T(i(p+1) + k, ip + k)$ with $\text{str}(Y) = t$, $0 \le t < n$. There can be at most one thread of any given order satisfying $0 \le \text{str}(T) < n$, and so $S_i$ must be the thread $Y$. In summary:

$$\text{ord}(S_i) = ip + k \qquad - (11)$$
$$\text{str}(S_i) = t \qquad - (12)$$
$$\text{len}(S_i) = (n+1) - i \qquad - (13)$$

Now we can complete the proof as four separate cases: $n$ and $s$, even or odd.

**When $n$ is even and $s$ is even**, we write $n = 2m$, $s = 2u$, and choose thread $S_m$:

(10) gives: $ms = kn + t \Rightarrow 2mu = 2mk + t$; so $k = u$, $t = 0$ and:

$\text{ord}(S_m) = mp + u$

$\text{str}(S_m) = 0$

$\text{len}(S_m) = (n + 1) - m$

$\text{str}(S_m) + \text{len}(S_m) = (n + 1) - m \le n + 1 = \text{len}(T_0)$; so $S_m$ is covered by $T_0$*.

So $2q < 2(\text{ord}(S_m)) = 2mp + 2u = pn + s = a_2$; so $q < a_2/2$.

But $n = C_1 < a_2/2$, so $n + q < a_2$ as required.

The thread $S_m$ appears as follows with respect to $T_0$

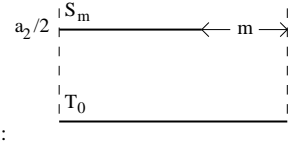

**When $n$ is even and $s$ is odd,** we write $n = 2m$, $s = 2u + 1$, and choose thread $S_m$:

(10) gives: $ms = kn + t \Rightarrow m(2u + 1) = 2mk + t$; so $k = u$, $t = m$ and:

$\quad \text{ord}(S_m) = mp + u$

$\quad \text{str}(S_m) = m$

$\quad \text{len}(S_m) = (n + 1) - m$

$\text{str}(S_m) + \text{len}(S_m) = n + 1 = \text{len}(T_0)$; so $S_m$ is covered by $T_0$*.

So $2q < 2(\text{ord}(S_m)) = 2mp + 2u = pn + s - 1 < a_2$; so $q < a_2/2$.

But $n = C_1 < a_2/2$, so $n + q < a_2$ as required.

The thread $S_m$ appears as follows with respect to $T_0$

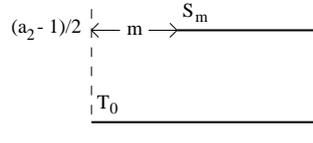

* Note that $m > 0$ (and so $S_m$ and $T_0$ are different threads) since when $n = C_1 = 0$, the stride generator is of order 0 by Lemma 15.

**When $n > 1$ is odd and $s$ is odd,** we write $n = 2m - 1$, $s = 2u + 1$, and choose thread $S_{m-1}$:

(10) gives: $(m-1)s = k(2m-1) + t \Rightarrow (2um - 2u + m - 1 - 2km + k) = t$. Substituting $k = u$ we get $t = m - u - 1$, and we now show that $0 \leq t < n$:

$\quad s < n \Rightarrow 2u+1 < 2m-1 \Rightarrow 2u < 2m-2 \Rightarrow u < m-1 \Rightarrow m-u-1 > 0 \Rightarrow t > 0$

$\quad u \geq 0 \Rightarrow t = m-u-1 \leq m-1 \leq (2m-2)/2 < (2m-1)/2 = n/2 < n$

So: $\text{ord}(S_{m-1}) = (m - 1)p + u$

$\quad \text{str}(S_{m-1}) = m - u - 1$

$\quad \text{len}(S_{m-1}) = n - m + 2$

$\text{str}(S_{m-1}) + \text{len}(S_{m-1}) = n - u + 1 \leq n + 1$; so $S_{m-1}$ is covered by $T_0$.

Since $n > 1$, $m > 1$ and so $S_{m-1}$ and $T_0$ are different threads; so $q < \text{ord}(S_{m-1})$. (If $m = 1$, $S_{m-1}$ is the same thread as $T_0$ and the cover argument is not applicable; this is why the case $n = 1$ must be dealt with specially.)

So $2q < 2(\text{ord}(S_{m-1})) = 2p(m-1) + 2u = p(2m-1) + 2u - p = pn + s - 1 - p = a_2 - (p + 1) < a_2$; so $q < a_2/2$.

But $n = C_1 < a_2/2$, so $n + q < a_2$ as required.

The thread $S_m$ appears as follows with respect to $T_0$

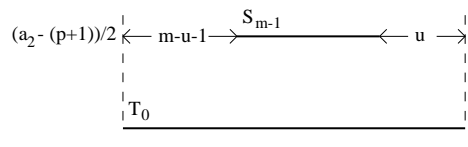

**When $n > 1$ is odd and $s$ is even,** we write $n = 2m - 1$, $s = 2u$, and choose thread $S_m$:

(10) gives: $2mu = k(2m - 1) + t$; so $k = u$, $t = u$ is the solution (since $0 \leq u \leq s < n$) and:

$\quad \text{ord}(S_m) = pm + u$

$$\text{str}(S_m) = u$$
$$\text{len}(S_m) = n + 1 - m = m$$

$\text{str}(S_m) + \text{len}(S_m) = u + m = s/2 + (n+1)/2 < n/2 + (n+1)/2 < n + 1$; so $S_m$ is covered by $T_0$; note that $S_m$ is a different thread from $T_0$ because $n = C_1 > 0 \Rightarrow m > 0$.

*For $p = 2$:*

$2q < 2(\text{ord}(S_m)) = 4m + 2u = 2n + s + 2 = a_2 + 2$; so $2q \leq a_2 + 1$.

By Lemma 16, $2n = 2C_1 < a_2 - 1$, so $2n + 2q < 2a_2 \Rightarrow n + q < a_2$ as required.

The thread $S_m$ appears as follows with respect to $T_0$

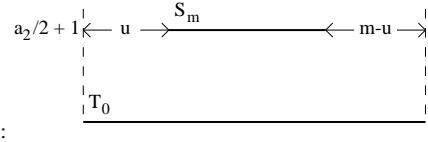

*For $p \geq 3$:*

$3q < 3(\text{ord}(S_m)) = 3pm + 3u = 2a_2 - p(m-2) - u \leq 2a_2$ provided that $m \geq 2$; once again the case $n = 1$ must be dealt with separately. So if $m > 1$, $q < (2/3)a_2$.

But $n = C_1 < a_2/p \leq a_2/3$; so $n + q < a_2$ as required.

The thread $S_m$ appears as follows with respect to $T_0$

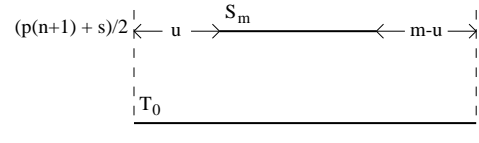

### *When $n = 1$:*

In this case, the stride generator is always canonical; for if it were not, (1) and (3) above would both be satisfied, which leads to a contradiction:

$(1) \Rightarrow (p+1)C_1 > a_2 \Rightarrow p > a_2 - 1$
$(3) \Rightarrow p(C_1) < a_2 - 1 \Rightarrow p < a_2 - 1$

This completes the proof for the case $C_1 < a_2/2$, $C_2 = 1$.

### *2.4 The ascending staircase: $C_1 > a_2/2$*

We know from Lemma 14 that any non-canonical fundamental stride generator with $C_1 > a_2/2$ has the following form of thread diagram for some $p \geq 2$:

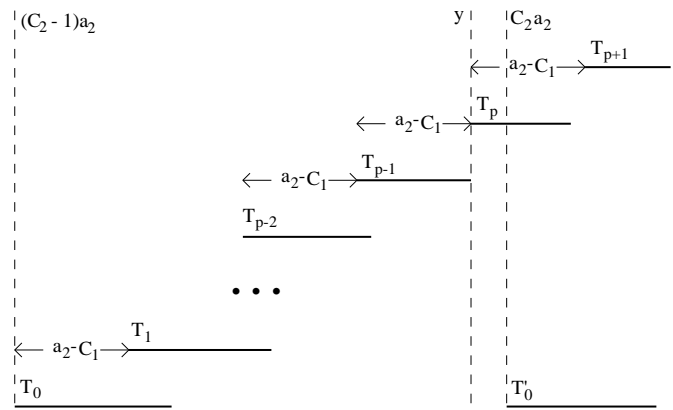

Ascending staircase: case (1)

We know that $T_{p+1}$ crosses the end of $T'_0$ because:

$str(T_{p+1}) > end(T_p) \geq str(T'_0)$, and $str(T_{p+1}) \leq end(T'_0)$ because:

$end(T'_0) = str(T'_0) + len(T'_0) - 1 = str(T'_0) + len(T_0) - 2 > str(T'_0) + (a_2 - C_1) - 2$

$str(T_{p+1}) = str(T_p) + (a_2 - C_1) < str(T'_0) + (a_2 - C_1)$

So $end(T'_0) \geq str(T'_0) + (a_2 - C_1) - 1$, and $str(T_{p+1}) \leq str(T'_0) + (a_2 - C_1) - 1$.

So $T_{p+1}$ is covered by $T'_0$ (and hence the stride generator is canonical) unless $end(T_{p+1}) > end(T'_0)$.

### 2.4.1 General bounds

We obtain bounds for $(a_2 - C_1)$ and a formula for n as follows:

We see immediately that $(p+1)(a_2 - C_1) > a_2 > p(a_2 - C_1)$. Here we need only the upper bound for $(a_2 - C_1)$, because we develop a better lower bound below in (3):

$$p(a_2 - C_1) < a_2 \qquad - (1)$$

Next, we note $len(T_{i+1}) = len(T_i) - C_2$ for $i \geq 0$, so:

$len(T_{p-1}) = len(T_0) - (p-1)C_2 = n - C_2 + 2 - (p-1)C_2$; but $len(T_{p-1}) = a_2 - C_1$, so:

$$n = (a_2 - C_1) + pC_2 - 2 \qquad - (2)$$

Since $T_{p+1}$ crosses the end of $T'_0$, we have $end(T_{p+1}) > end(T'_0)$:

$end(T'_0) = C_2 a_2 + n - C_2 = C_2 a_2 + (a_2 - C_1) + (p-1)C_2 - 2$

$end(T_{p+1}) = (C_2 - 1)a_2 + (p+2)(a_2 - C_1) - 2C_2 - 1$

So: $end(T_{p+1}) > end(T'_0) \Rightarrow (p+1)(a_2 - C_1) > a_2 + (p+1)C_2 - 1$, or:

$$(p + 1)((a_2 - C_1) - C_2) > a_2 - 1 \qquad - (3)$$

We can now derive an upper bound for $C_2$ as follows:

From (3), $(p+1)C_2 < (p+1)(a_2 - C_1) - a_2 + 1$

From (1), $p(a_2 - C_1) < a_2 \Rightarrow (p+1)(a_2 - C_1) < ((p+1)/p)a_2$

So: $(p+1)C_2 < ((p+1)/p)a_2 - a_2 + 1 = a_2/p + 1$

$$\Rightarrow \qquad C_2 < a_2/(p(p+1)) + 1/(p+1) \qquad - (4)$$

### 2.4.2 The case for $C_2 \geq 2$

In this section, we assume $C_2 \geq 2$.

We obtain an upper bound for q as follows:

Let $T_j$ be the thread of order j such that $(C_2 - 1)a_2 \leq str(T_j) < C_2 a_2$:

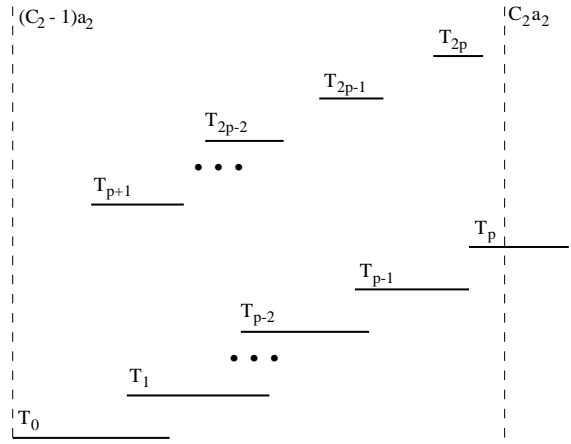

If $T_{j+1}$ and $T_j$ are two steps on the same staircase, $len(T_{j+1}) = len(T_j) - C_2$. If they are at opposite ends of the range, $len(T_{j+1}) = len(T_j) - (C_2 - 1)$; this happens at most every pth thread, starting with $T_{p+1}$. So:

$$len(T_i) \le len(T_0) - iC_2 + (i-1)/p \le len(T_0) - iC_2 + (i-1)/2 \text{ since } p \ge 2$$
$$\Rightarrow len(T_i) \le n - C_2 + 2 - iC_2 + (i-1)/2 = n - C_2 + 3/2 - i(C_2 - 1/2)$$

So $len(T_i) \le 0$ if $i(2C_2 - 1) \ge 2n - 2C_2 + 3$; and so $q < (2n - 2C_2 + 3)/(2C_2 - 1)$, or:

$$q < (2n + 2)/(2C_2 - 1) - 1 \qquad - (5)$$

We can now derive an upper bound for $n + q$:

From (5): $n + q < n + 2n/(2C_2 - 1) - (2C_2 - 3)/(2C_2 - 1)$
$$= n((2C_2 + 1)/(2C_2 - 1)) - (2C_2 - 3)/(2C_2 - 1)$$

From (1): $(a_2 - C_1) < a_2/p$; substituting in (2) we have: $n < a_2/p + pC_2 - 2$

So: $\qquad n + q < (a_2/p + pC_2 - 2)((2C_2 + 1)/(2C_2 - 1)) - (2C_2 - 3)/(2C_2 - 1) \qquad - (6)$

From (4): $pC_2 < a_2/(p+1) + p/(p+1) < a_2/(p+1) + 1$; so:
$$n + q < (a_2/p + a_2/(p+1) - 1)((2C_2 + 1)/(2C_2 - 1)) - (2C_2 - 3)/(2C_2 - 1), \text{ which gives:}$$
$$n + q < (a_2/p + a_2/(p+1))((2C_2 + 1)/(2C_2 - 1)) \qquad - (7)$$

We are now ready to prove that $n + q \le a_2$; we take the cases $p \ge 3$ and $p = 2$ separately.

*When $p \ge 3$:*

From (7) we have:
$$n + q < (5/3)((2p+1)/(p(p+1)))a_2 \text{ since } C_2 \ge 2$$
$$\Rightarrow n + q < (35/36)a_2 < a_2 \qquad \text{for all } p \ge 3$$

*When $p = 2$:*

From (7) we have:
$$n + q < (a_2/2 + a_2/3)((2C_2 + 1)/(2C_2 - 1))$$
$$= (5/6)a_2((2C_2 + 1)/(2C_2 - 1)) \le (65/66)a_2 < a_2 \text{ for } C_2 \ge 6$$

This leaves $C_2 = 2, 3, 4$ and $5$ to be considered; we note that, from (4):
$$2C_2 < (a_2 + 2)/3 \Rightarrow a_2 > 6C_2 - 2 \qquad - (8)$$

For $C_2 = 5$, we substitute directly in (6):

$$n + q < (a_2/2 + 8)(11/9) - (7/9) = (11a_2 + 162)/18, \text{ which is } \leq a_2 \text{ so long as } a_2 \geq 24.$$

From (8) we find $a_2 > 28$, which is sufficient.

For $C_2 = 4$, we similarly substitute directly in (6):

$$n + q < (a_2/2 + 6)(9/7) - (5/7) = (9a_2 + 98)/14, \text{ which is } \leq a_2 \text{ so long as } a_2 \geq 20.$$

From (8) we find $a_2 > 22$, which is sufficient.

For $C_2 = 3$, we again substitute directly in (6):

$$n + q < (a_2/2 + 4)(7/5) - (3/5) = (7a_2 + 50)/10, \text{ which is } \leq a_2 \text{ so long as } a_2 \geq 17.$$

From (8) we find $a_2 > 16$, which is just sufficient.

For $C_2 = 2$, (6) gives:

$$n + q < (a_2/2 + 2)(5/3) - (1/3) = (5a_2 + 18)/6, \text{ which is } \leq a_2 \text{ so long as } a_2 \geq 18.$$

But (8) requires only that $a_2 > 10$, so we have $11 \leq a_2 \leq 17$ to consider.

From (5) we have $q < (2n - 1)/3 \Rightarrow n + q < (5n - 1)/3$; and (2) gives $n = (a_2 - C_1) + 2$.
Substituting, we have $n + q < (5a_2 - 5C_1 + 9)/3$, which is $\leq a_2$ so long as $C_1 \geq (2a_2 + 9)/5$.

We now construct the following table where X is the smallest value $C_1 > a_2/2$, and Y is the smallest value $C_1 \geq (2a_2 + 9)/5$:

| $a_2$: | 11 | 12 | 13 | 14 | 15 | 16 | 17 |
|---|---|---|---|---|---|---|---|
| X: | 6 | 7 | 7 | 8 | 8 | 9 | 9 |
| Y: | 7 | 7 | 7 | 8 | 8 | 9 | 9 |

We have only to consider cases where $X \leq C_1 < Y$, which leaves the single case $p = 2$, $C_2 = 2$, $a_2 = 11$, $C_1 = 6$ which results in $n = 7$ and $q < 13/3$, or $q \leq 4$; so $n + q \leq 11 = a_2$ as required.

(This is the second case where we show only that $n + q \leq a_2$. As shown above, the only doubtful case is $A = \{1, 11, 28\}$ which turns out to be a canonical stride generator.)

This completes the proof for the case $C_1 > a_2/2$, $C_2 \geq 2$.

### 2.4.3 The case for $C_2 = 1$

With $C_1 > a_2/2$, we find that $(a_2 - C_1)$ plays a similar role to that of $C_1$ when $C_1 < a_2/2$; so we write:

$$n' = (a_2 - C_1) \qquad - (8)$$

From (2), we have:

$$n = n' + p - 2 \qquad - (9)$$

and, since $p(a_2 - C_1) < a_2 < (p+1)(a_2 - C_1)$, we can write:

$$a_2 = pn' + s \text{ where } 1 \leq s < n' \qquad - (10)$$

As before, let $T_j$ be the thread of order j such that $0 \leq \text{str}(T_j) < a_2$. We now consider those threads $S_i$ which satisfy $0 \leq \text{str}(S_i) < n'$:

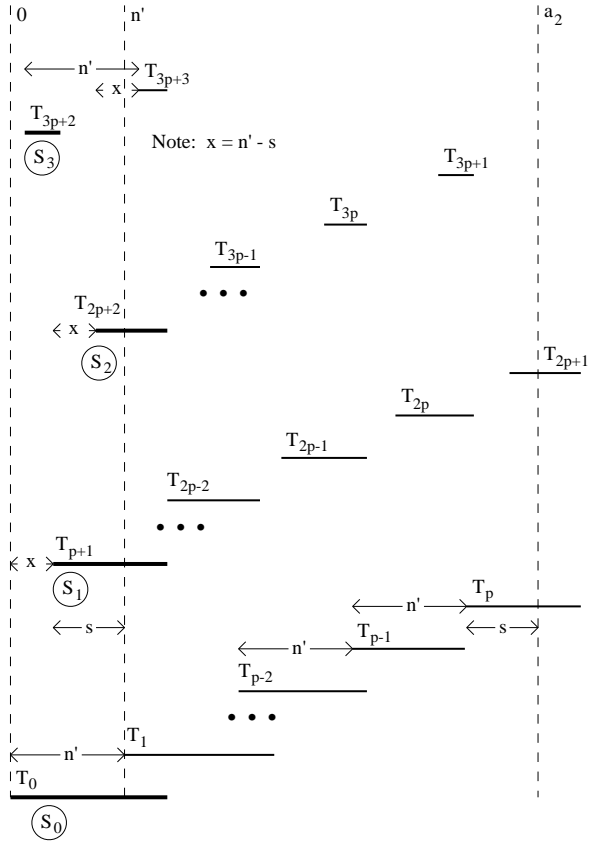

Since $\text{len}(T_0) = n+1 = n' + p - 1 \leq n'+1$, each thread $S_i$ at least satisfies $\text{str}(T_0) \leq \text{str}(S_i) < \text{end}(T_0)$. We will show below that there is always a thread $S_i$ that satisfies $\text{end}(S_i) \leq \text{end}(T_0)$ and so is covered by $T_0$, thus providing an upper bound for $q$; it turns out that this bound is sufficient to show that $n + q < a_2$.

First we derive formulae for the order, position and length of thread $S_i$, by observing that $S_i$ is derived from the thread $X = T(2(ip+1)-i, ip+1)$ as follows:

We have: $\text{str}(X) = (2ip + 2 - i)a_2 - (ip - 1)a_3$
$= (2ip + 2 - i)a_2 - (ip + 1)(2a_2 - n')$      since $a_3 = a_2 + C_1 = 2a_2 - n'$
$= (ip + 1)n' - ia_2 = (ip + 1)n' - i(pn' + s) = n' - is$

and: $\text{len}(X) = n + ip+1 - 2(ip+1) + i + 1 = n - ip + i = n - i(p - 1)$

With $C_2 = 1$, a thread $X_1 = T(c, i)$ implies the existence of a further thread $X_2 = T(c+1, i+1)$ of the same length where $\text{str}(X_2) = \text{str}(X_1) - C_1$; and from $X_2$ we can derive $X_3 = T(c+2, i+1)$ of length one less with $\text{str}(X_3) = \text{str}(X_2) + a_2 = \text{str}(X_1) + (a_2 - C_1) = \text{str}(X_1) + n'$. Let us write:

$$is = kn' - t \text{ where } 0 \leq t < n' \quad - (11)$$

Then from X we derive thread $Y = T(2(ip+k)-i, ip+k)$ with $\text{str}(Y) = t$, $0 \leq t < n'$. There can be at most one thread of any given order satisfying $0 \leq \text{str}(T) < n'$, and so $S_i$ - if it exists - must be the thread Y. In summary:

$$\text{ord}(S_i) = ip + k \quad - (12)$$
$$\text{str}(S_i) = t \quad - (13)$$
$$\text{len}(S_i) = n - i(p - 1) - k + 1 \quad - (14)$$

Once we have found a thread $S_m$ with $m > 0$ that is covered by $T_0$, it is sufficient to show that $n + \text{ord}(S_m) - 1 < a_2$ in order to prove that $n + q < a_2$. Using (12), (9) and (10), this is equivalent to showing that $n' + p - 2 + mp + k - 1 < pn' + s$, or:

$$pn' + s - n' - p + 3 - mp - k > 0 \qquad - (15)$$

Now we can complete the proof as four separate cases: n' and s, even or odd.

**When n' is even and s is even,** we write $n' = 2m$, $s = 2u$, and choose thread $S_m$:

  (11) gives: $ms = kn' - t \Rightarrow 2mu = 2mk - t$; so $k = u$, $t = 0$ and:

  $\quad \text{ord}(S_m) = mp + u$

  $\quad \text{str}(S_m) = 0$

  $\quad \text{len}(S_m) = n - m(p - 1) - u + 1$

  $\text{str}(S_m) + \text{len}(S_m) = (n + 1) - m(p - 1) - u \leq n + 1 = \text{len}(T_0)$; so $S_m$ is covered by $T_0$*.

  (15) gives: $2mp + 2u - 2m - p + 3 - mp - u = m(p - 2) - p + u + 3$

  $\quad \geq p - 2 - p + u + 3 = u + 1 > 0$ since $m \geq 1$ and $p \geq 2$.

**When n' is even and s is odd,** we write $n' = 2m$, $s = 2u - 1$, and choose thread $S_m$:

  (11) gives: $ms = kn' - t \Rightarrow m(2u - 1) = 2mk - t$; so $k = u$, $t = m$ and:

  $\quad \text{ord}(S_m) = mp + u$

  $\quad \text{str}(S_m) = m$

  $\quad \text{len}(S_m) = n - m(p - 1) - u + 1$

  $\text{str}(S_m) + \text{len}(S_m) = n - mp + 2m - u + 1 = n + 1 - m(p - 2) - u < n + 1 = \text{len}(T_0)$; so $S_m$ is covered by $T_0$*.

  (15) gives: $2mp + 2u - 1 - 2m - p + 3 - mp - u = m(p - 2) - p + u + 2$

  $\quad \geq p - 2 - p + u + 2 = u > 0$ since $m \geq 1$, $u \geq 1$ and $p \geq 2$.

  \* Note that $S_m$ and $T_0$ are different threads because $C_1 < a_2 \Rightarrow n' > 0 \Rightarrow m > 0$.

**When n' > 1 is odd and s is odd,** we write $n' = 2m + 1$, $s = 2u - 1$, and choose thread $S_m$:

  (11) gives: $m(2u-1) = k(2m+1) - t \Rightarrow (2km + k - 2um + m) = t$. Substituting $k = u$ we get $t = m + u$, and we now show that $0 \leq t < n'$:

  $\quad 0 \leq 2t = 2m + 2u = n' - 1 + s + 1 = n' + s < 2n'$

  So: $\text{ord}(S_m) = mp + u$

  $\quad \text{str}(S_m) = m + u$

  $\quad \text{len}(S_m) = n - m(p - 1) - u + 1$

  $\text{str}(S_m) + \text{len}(S_m) = m + n - m(p - 1) + 1 = n + 1 - m(p - 2) \leq n + 1$; so $S_m$ is covered by $T_0$.

  Since $n' > 1$, $m > 0$ and so $S_m$ and $T_0$ are different threads; so $q < \text{ord}(S_m)$. (If $m = 0$, $S_m$ is the same thread as $T_0$ and the cover argument is not applicable; this is why the case $n' = 1$ must be dealt with specially.)

  (15) gives: $p(2m + 1) + 2u - 1 - 2m - 1 - p + 3 - mp - u = m(p - 2) + u + 1$

  $\quad \geq p - 2 + u + 1 \geq u + 1 > 0$ since $m \geq 1$ and $p \geq 2$.

**When n' > 1 is odd and s is even,** we write $n' = 2m + 1$, $s = 2u$, and choose thread $S_m$:

  (11) gives: $2mu = k(2m + 1) - t$; so $k = u$, $t = u$ is the solution and:

  $\quad \text{ord}(S_m) = pm + u$

$$\text{str}(S_m) = u$$
$$\text{len}(S_m) = n - m(p - 1) - u + 1$$

$\text{str}(S_m) + \text{len}(S_m) = n - m(p - 1) + 1 \leq n + 1$; so $S_m$ is covered by $T_0$; note that $S_m$ and $T_0$ are different threads because $n' > 1 \Rightarrow m > 0$.

(15) gives: $p(2m + 1) + 2u - 2m - 1 - p + 3 - mp - u = m(p - 2) + u + 2$

$\geq p - 2 + u + 2 = u + p > 0$ since $m \geq 1$ and $p \geq 2$.

*When $n' = 1$:*

$n' = 1 \Rightarrow C_1 = a_2 - 1$, and so the stride generator is of order 0 by Lemma 15; this contradicts our assumption that $p \geq 2$, and so this case cannot arise.

This completes the proof for the case $C_1 < a_2/2$, $C_2 = 1$.

**Acknowledgment**

I would like to express my thanks to Professor Ernst Selmer for persevering with [1] and observing that $h_1 \leq h_0$ implies and is implied by the conjecture that every underlying stride generator is canonical. He then encouraged me to see if a proof of my conjecture could be developed with the aid of thread diagrams - and this is the result.

**References**

[1] Challis, M.F., "The Postage Stamp Problem: Formulae and proof for the case of three denominations", Storey's Cottage, Whittlesford, Cambridge (1990).

[2] Hofmeister, G., "Asymptotische Abschätzungen für dreielementige Extremalbasen in natürlichen Zahlen", Journal für die reine und angewandte Mathematik 232 (1968), pp.77-101.

[3] Rodseth, O.J., "On h-bases for n", Math. Scand. 48 (1981), pp.165-183.

[4] Selmer, E.S., "The Local Postage Stamp Problem Part 1: General Theory", Research monograph No. 42, Department of Mathematics, University of Bergen, Norway (1986).

[5] Windecker, R., "Zum Reichweitenproblem", Dissertation, Math. Inst., Joh. Gutenberg-Univ., Mainz (1978).

# Appendix A   Historical information

On March 14th 1995 I received a letter from Selmer (after he had read [1]) in which he suggested that it might be possible to show $h_1 \leq h_0$ for $k = 3$ using my concepts of threads and stride generators. He observed that $h_1 \leq h_0$ is equivalent to Conjecture 1 in [1]: that is, $h_1 \leq h_0$ iff "all underlying stride generators are canonical".

I read his notes on $h_1$ and $h_2$, and also re-read Chapter 2 of [1]. From Conjecture 1 and the proof of Theorem 236, we know that any non-canonical underlying stride generator must have $n + q > a_2$, so it is sufficient to show that for any non-canonical stride generator $n + q \leq a_2$: this is, of course, the main theme of this paper.

I also became interested in showing that $h_2 \leq h_0$, and I was soon able to show that this followed immediately if I could show that in any canonical stride generator every thread of order $\geq p+1$ is covered by one of order $\leq p$. My first proof of this - completed in August 1995 - was based on Theorems 220 and 225 in [1], and is reproduced in version 0.02 of this document (which was sent to Selmer on January 26th 1996). The more elegant proof given in version 0.03 was developed later during February 1996.

But the real challenge was to show that $n + q \leq a_2$ for any non-canonical stride generator, and I started in September 1995 by developing a program (now called NCSTRIDES) which systematically generates all stride generators satisfying $a_3 < a_2^2$. Note that this is sufficient to include all interesting ones, since we know from [1], Theorem 214, that A is of order 0 when $C_1 \geq a_2 - C_2$; and so all stride generators with $a_3 \geq a_2^2$ are of order 0 and hence canonical ([1], Theorem 217). The original program was easily modified to list only non-canonical stride generators, and after a few hours' computing I had a file containing details of all those satisfying $a_2 \leq 138$; there are 74541 of them altogether. I also wrote an auxiliary program, PROCSTRDS, to read in and process these details; this was easily modified as required to "filter" the input so that I could investigate various hypotheses about the properties of non-canonical stride generators.

My initial experiments showed that:
- $n + q < a_2$         (as expected!)
- $q \leq 2p$         (with $q = 2p$ quite common)

This suggested that it might be possible to prove the stronger result $n + 2p \leq a_2$, but I soon found that:

$n + 2p \geq a_2$    for some *canonical* stride generators: $\{1, 8, 49\} = SG(11, 1)$, $\{1, 8, 9\} = SG(1, 6)$, $\{1, 14, 39\} = SG(10, 4)$

(See also [1], Conjecture 1, where we show that $n - a_2$ can become arbitrarily large for order zero stride generators)

$n + 2p \geq a_2$    for some *non-canonical* stride generators: $\{1, 65, 98\} = SG(19, 28)$ has $q = 30$

But, intriguingly, it *does* seem to be the case that $n + 2p < a_2$ for all non-canonical stride generators with $C_2 \geq 2$ (although I have not tried to prove this).

I noted next that:
- q and p are often related arithmetically in a simple way - but not always: $\{1, 93, 104\} = SG(6, 24)$ has $q = 41$
- $q \geq 4$    (this follows because a stride generator is canonical if $p = 0$ or $1$, and so $q > p+1 \Rightarrow q > 3$)
- For a given value of $a_2$, there is a maximum value of $C_2$ for which non-canonical stride generators exist; this increases as $a_2$ increases, but it looks as if $\max(C_2) \sim a_2/6$:

| $a_2$ | $\max(a_3)$ | $\max(C_2)$ | $C_2/a_2$ |
|---|---|---|---|
| 92 | 1411 | 15.34 | 0.166... |
| 138 | 3106 | 22.51 | 0.163... |

(I later showed that $C_2 < a_2/(p(p+1)) + 1 < a_2/6 + 1$ since $p \geq 2$; see sections 2.3.1 and 2.4.1 of this paper.)

I clearly needed to know more about the behaviour of non-canonical stride generators, and I looked again at the "series" of stride generators $SG(A, n_i, p_i)$ first described in [1], Theorem 232. This led me to the idea of the *fundamental* stride generator which has the form of either an ascending or descending staircase, and I used the program EXP3 (developed for [1]) to print out thread diagrams for each stride generator $\{1, 30, a_3\}$ for $31 \leq a_3 \leq 150$ to examine this hypothesis in more detail.

Since the fundamental stride generator has maximum n, any break in the fundamental or in any stride generator derived from it (ie in the same series) must be of order less than or equal to that ($q_{max}$) of the smallest thread in the fundamental. This idea turned out to be sufficient to prove $n + q < a_2$ for $C_2 \geq 2$, although the final special cases took some time to pin down. As shown in this paper, the cases $C_1 < a_2/2$ (descending staircase) and $C_1 > a_2/2$ (ascending staircase) are treated separately.

The case $C_2 = 1$ would not yield to this simple approach (indeed, $q_{max} = a_2$ for all $C_1 < a_2/2$ when $C_2 = 1$), and further investigation was needed. I first concentrated on $C_1 < a_2/2$ - correctly expecting this to be the most difficult case - and hit on the idea of looking just at those threads $S_1, S_2$ ... which overlapped $T_0$. If thread $S_i$ is covered by $T_0$, then $q_{max} \leq ord(S_i) - 1$; and so I started looking for classes of thread $S_i$ with this property which could also be shown to satisfy $n + ord(S_i) - 1 \leq a_2$. A new program - INVTHRD - was developed to list details of these threads for selected stride generators, and examination of the output immediately showed certain tantalising patterns; but it was some time before I was able to interpret these fully:

For $C_2 = 1$, $C_1 < a_2/2$  we have  $n = C_1$  and write  $a_2 = pn + s$; then:
- There are n threads  $S_1, S_2, ... S_n$  which satisfy  $0 \leq str(S_i) < n$
- $ord(S_n) = a_2$,  $str(S_n) = 0$,  $len(S_n) = 1$
- These threads are spaced vertically as equally as possible, so that  $ord(S_{i+1}) - ord(S_i) = p$  or $p+1$.

The threads $S_i$ fall naturally into groups of roughly equal size, each with a thread with a (locally) minimum offset $str(S_i)$; we denote these *minimal* threads $M_j$.
- If $s < n/2$, there are s groups, each one consisting of threads $S_i$ where $str(S_i) = str(S_{i-1}) + s$
- If $s > n/2$, there are (n-s) groups, each one consisting of threads $S_i$ where $str(S_i) = str(S_{i-1}) - (n-s)$

For example:

A = $\{1, 30, 37\} = SG(7, 4)$;  s = 2   (s < n/2)

| i | $ord(S_i)$ | $str(S_i)$ | |
|---|---|---|---|
| 1 | 4 | 2 | |
| 2 | 8 | 4 | |
| 3 | 12 | 6 | |
| 4 | 17 | 1 | |
| 5 | 21 | 3 | |
| 6 | 25 | 5 | |
| 7 | 30 | 0 | |

The first group is $\{S_1, S_2, S_3, S_4\}$  with  $M_1 = S_4$;  the second group is $\{S_5, S_6, S_7\}$  with  $M_2 = S_7$.

and:

A = $\{1, 31, 43\} = SG(12, 2)$; s = 7 => (n-s) = 5  (s > n/2)

| i | $ord(S_i)$ | $str(S_i)$ | |
|---|---|---|---|
| 1 | 2 | 7 | |
| 2 | 5 | 2 | $M_1$ |
| 3 | 7 | 9 | |
| 4 | 10 | 4 | $M_2$ |
| 5 | 12 | 11 | |
| 6 | 15 | 6 | |
| 7 | 18 | 1 | $M_3$ |
| 8 | 20 | 8 | |
| 9 | 23 | 3 | $M_4$ |
| 10 | 25 | 10 | |
| 11 | 28 | 5 | |
| 12 | 31 | 0 | $M_5$ |

Later, I found that threads $S_i$ for $C_1 > a_2/2$ have similar properties:

For $C_2 = 1$, $C_1 > a_2/2$ we write $n' = a_2 - C_1$, $a_2 = pn' + s$; then $n' = n - p + 2$ and:
- If $s < n'/2$, each group consists of threads $S_i$ where $str(S_i) = str(S_{i-1}) - s$
- If $s > n'/2$, each group consists of threads $S_i$ where $str(S_i) = str(S_{i-1}) + (n-s)$

The threads $S_i$ are spaced vertically as equally as possible, with $ord(S_{i+1}) - ord(S_i) = p$ or $p+1$, but because their length decreases more rapidly as i increases than when $C_1 < a_2/2$, the last thread $S_i$ has order $< a_2$.

What are the chances that we can always find a minimal thread $M_i$ that is covered by $T_0$ and whose order is sufficiently small that $n + ord(M_i) - 1 \leq a_2$?
As i increases so $len(M_i)$ decreases, making it more likely that $M_i$ is covered; so we try to choose $ord(M_i)$ as large as possible.
Note that $str(M_i) < s$ (or n-s as appropriate) for all i, so our chances improve when s (or n-s) is small relative to n/2.
On the other hand, we have a greater choice of minimal threads when the number of groups is large - that is, as s (or n-s) approaches n/2.
It turns out that all is well for "reasonable" values of s, but extreme cases - notably $s = 1$ (or n-s = 1) where there is only one minimal thread $M_1 = S_n$ - have to be dealt with specially.

The overall approach to the $C_2 = 1$, $C_1 < a_2/2$ case can be summarised as follows:
- Choose a suitable thread $S_i$
- Show that $S_i$ is covered by $T_0$
- Show that $n + ord(S_i) - 1 \leq a_2$

The question remains as to how to choose the thread $S_i$. Possibilities I investigated included:
a) Choose the thread with $str(S_i) = 1$.
   This does not work, because although it is always covered by $T_0$, its order is sometimes too great.
b) Choose the first minimal thread $M_1$.
   This fails because $M_1$ is not always covered by $T_0$.
c) Choose the minimal thread $M_j$ with highest order $\leq Q$ where $Q = a_2/2 + 1$ for $p = 2$, and $Q = 2a_2/3$ for $p \geq 3$.
   (The reason for these choices of Q is given below)
d) Choose the thread $S_{n/2}$ (or thereabouts).

Approach (c) was my first success. The proof is complex, involving separate cases according as $C_1 < a_2/2$, $C_1 > a_2/2$; $p = 2$, $p \geq 3$; and $s < n/2$, $s > n/2$.
Once I was satisfied that (c) could be made to work, I wrote to Selmer (November 29th 1995) with an outline of my proof, and then proceeded to sort out the details. These proved trickier than expected, and it was then that approach (d) occurred to me; this is reproduced in the main body of this paper, and was sent to Selmer on January 26th 1996. Note that (d) is much simpler than (c) because there is no requirement for the chosen thread to be *minimal*. In February 1996 I returned to the details of (c), and - although of historical interest only - these are reproduced in Appendix B below.

One approach that I followed for $C_2 = 1$ (and which is used in many of the sub-cases for both approaches (c) and (d)) is to "divide and conquer" by determining separate bounds for n and q which, when taken together, show that $n + q \leq a_2$; for example, for $C_1 < a_2/2$ we have $n = C_1 < a_2/2$, and so if we can show that $q < a_2/2$, we are home and dry.
Experiments using PROCSTRDS suggested the following to be true:
  For all non-canonical stride generators:
  - $q < 2a_2/3$
  - $q > a_2/2$ only when $C_2 = 1$, $C_1 < a_2/3$ (ie $p \geq 3$)
      For example, $\{1, 3t+2, 3t+5\} = SG(3, n)$ has $q = 2n$; so as n tends to infinity, $q \to 2a_2/3$
  - The maximum value of $n + q$ seems to be around $2a_2/3$ as $a_2$ becomes large, but the worst case is:
      $\{1, 11, 14\} = SG(3, 3)$ with $q = 6$, where $n + q = 9/11 = 0.818...$

These observations suggest splitting the case $C_1 < a_2/2$ into two as follows:
  (i)  $0 < n = C_1 < a_2/3 \iff p \geq 3$; we have $n < a_2/3$, and must prove $q \leq 2a_2/3$.
  (ii) $a_2/3 < n = C_1 < a_2/2 \iff p \geq 2$; we have $n < a_2/2$, and must prove $q \leq a_2/2$.
This explains the choices for Q above.

Approaches (c) and (d) identify an upper bound for q as one less than the order of a thread covered by $T_0$, and so it seemed sensible to check out the properties of $q_{max}$ - one less than the order of the first such thread. It turns out that $q_{max}$ is a sharp bound for q, and for non-canonical stride generators I found:
- When $C_1 \geq a_2/3$, we find $q_{max} \leq a_2/2$ with equality only when $C_1 = a_2/2 - 1$ or $C_1 = (a_2+1)/3$. Under some conditions, $q = q_{max}$:

    | $a_2$ | $C_1$ | $q = q_{max}$ | |
    |---|---|---|---|
    | 4t | $a_2/2 - 1$ | $a_2/2$ | - (A1) |
    | 6t+2 | $(a_2+1)/3$ | $a_2/2$ | - (A2) |

    and so both conditions arise when $a_2 = 12t+8$, as:
        $A = \{1, 12t+8, 18t+11\} = SG(3t+3, 6t+2)$, $q = 6t+4$
        $A = \{1, 12t+8, 16t+11\} = SG(2t+3, 6t+1)$, $q = 6t+4$

- When $C_1 \leq a_2/3$, we find $q_{max} \leq (2a_2 - 4)/3$ with equality only when $C_1 = 3$:
        $A = \{1, 3t+2, 3t+5\} = SG(3, t)$, $q_{max} = q = 2t$  (cf above) - (B)

It is interesting to check that these observations are consistent with the results of section 2.3.3:
  (A1) gives $n = a_2/2 - 1 \implies a_2 = 2n+2 \implies s = 2$
  (A2) gives $n = (a_2+1)/3 \implies a_2 = 2n + (n-1) \implies s = n-1$
  Only the case $n > 1$, n odd and s even allows the possibility that $q = a_2/2$, and when $p = 2$ we find $ord(S_m) = a_2/2 + 1$ - which is just consistent with $q_{max} \leq a_2/2$. So both (A1) and (A2) also require n to be odd.
  (B) gives $n = 3$, and only the case n odd, s even, $p \geq 3$ allows $q_{max} > a_2/2$. The highest value for $ord(S_m)$ arises when $u = 1$ and we have:

$n = 3 \Rightarrow m = 2$; $u = 1 \Rightarrow s = 2$; so $a_2 = pn + s = 3p + 2 \Rightarrow p = 3$.
So $\text{ord}(S_m) = (2a_2 - 1)/3$, which is just consistent with $q_{max} \leq (2a_2 - 4)/3$.

Further experiments were undertaken for the case $C_2 = 1$, $C_1 > a_2/2$, where we already know that $q < a_2/2$; we found:
For all non-canonical stride generators with $C_2 = 1$, $C_1 > a_2/2$:
- $q \leq (a_2 - 1)/2$ with equality only for: $\{1, 4t+1, 6t+2\} = SG(t+2, 2t-2)$, $q = 2t$.
- $n \leq (a_2 - 1)/2$ with equality only for: $\{1, 2t+1, 3t+2\} = SG(t, 2)$, $q = 4$.

This result strongly suggests a proof split along the lines $n < a_2/2$ and $q < a_2/2$, and I soon managed to prove the former (see details in Appendix B). However, a demonstration that $q < a_2/2$ has proved more elusive, and the proof for $C_1 > a_2/2$ given in Appendix B is split into two parts as follows:

For $p = 2$, we show $q \leq a_2/2$ using techniques similar to those used in the $C_1 < a_2/2$ proof.

For $p \geq 3$, I have been unable to prove that $q \leq a_2/2$: I can only manage $q \leq a_2/2$ for $p \geq 4$, and $q < a_2/2 + (5/2)$ for $p = 3$. Instead I use a separate argument developed in February 1996 and derived from section 2.4.2 above.

# Appendix B   Alternative proof for $C_2 = 1$

The proof is structured as follows:

| | | | | | |
|---|---|---|---|---|---|
| 1 | $C_1 < a_2/2$ | | | | For $C_1 < a_2/2$: |
| 1.1 | | s = 1 | | | |
| 1.2 | | s = n/2 | | | |
| 1.3 | | s < n/2 | | | For s < n/2: |
| 1.3.1 | | | k > 2 | | $a_3 = a_2 + C_1$ |
| 1.3.1.1 | | | | p = 2 | $n = C_1$ |
| 1.3.1.1.1 | | | | | $s \leq k+2$    $a_2 = pn + s$    $2 \leq s < n$ |
| 1.3.1.1.2 | | | | | $s > k+2$    $n = ks + t$    $0 \leq t < s$ |
| 1.3.1.2 | | | | $p \geq 3$ | |
| 1.3.1.2.1 | | | | | $s \leq k+2$ |
| 1.3.1.2.2 | | | | | $s > k+2$ |
| 1.3.2 | | | k = 2 | | |
| 1.3.2.1 | | | | p = 2 | |
| 1.3.2.1.1 | | | | | s even |
| 1.3.2.1.2 | | | | | s odd |
| 1.3.2.2 | | | | $p \geq 3$ | |
| 1.4 | | s > n/2 | | | For s > n/2, we write s' = n-s: |
| 1.4.1 | | | s' = 1 | | $a_3 = a_2 + C_1$ |
| 1.4.1.1 | | | | n even | $n = C_1$ |
| 1.4.1.2 | | | | n odd, $n \geq 3$ | $a_2 = pn + s = (p+1)n - s'$    $1 \leq s < n$ |
| 1.4.1.2.1 | | | | | p = 2    $n = ks' + t$    $0 \leq t < s'$ |
| 1.4.1.2.2 | | | | | $p \geq 3$ |
| 1.4.1.3 | | | | n = 1 | |
| 1.4.2 | | | $s' \geq 2$ | | |
| 1.4.2.1 | | | | k > 2 | |
| 1.4.2.1.1 | | | | | $s' \leq k+1$ |
| 1.4.2.1.2 | | | | | $s' > k+1$ |
| 1.4.2.1.2.1 | | | | | p = 2 |
| 1.4.2.1.2.2 | | | | | $p \geq 3$ |
| 1.4.2.2 | | | | k = 2 | |
| 1.4.2.2.1 | | | | | p = 2 |
| 1.4.2.2.1.1 | | | | | $t \geq s'/3$ |
| 1.4.2.2.1.2 | | | | | $t < s'/3$ |
| 1.4.2.2.1.3 | | | | | Alternative method for all t |
| 1.4.2.2.2 | | | | | $p \geq 3$ |
| | | | | | |
| 2 | $C_1 > a_2/2$ | | | | For $C_1 > a_2/2$: |
| 2.1 | | p = 2 | | | |
| 2.1.1 | | | s = n'/2 | | $a_3 = a_2 + C_1$ |
| 2.1.2 | | | s < n'/2 | | $n' = a_2 - C_1$ |
| 2.1.2.1 | | | | $s \geq 18$ | $a_2 = pn' + s$    $0 \leq s < n'$ |
| 2.1.2.2 | | | | $2 \leq s \leq 17$ | $n' = ks + t$    $0 \leq t < s$ |
| 2.1.2.3 | | | | s = 1 | s' = n' - s |
| 2.1.3 | | | s > n'/2 | | |
| 2.1.3.1 | | | | $s' \geq 3$ | n = n' + p - 2 |
| 2.1.3.1.1 | | | | | k > 2 |
| 2.1.3.1.2 | | | | | k = 2 |
| 2.1.3.2 | | | | s' = 2 | |
| 2.1.3.3 | | | | s' = 1 | |
| 2.2 | | $p \geq 3$ | | | |
| 2.2.1 | | | $p \geq 5$ | | |
| 2.2.2 | | | p = 4 | | |
| 2.2.3 | | | p = 3 | | |

## *1  $C_1 < a_2/2$*

From section 2.3.3 we know that $n = C_1$, and that $pC_1 < a_2 < (p+1)C_1$; we write:

$a_3 = a_2 + C_1 = a_2 + n$
$a_2 = pn + s$    $1 \leq s < n$
$n = ks + t$    $0 \leq t < s$
$s' = n - s$

Clearly, $n < a_2/p \leq a_2/2$.

We use the notation $S_0 = T_0$. $S_1, S_2, ...$ to identify those threads which satisfy $0 \leq str(S_i) < n$ (see the diagram in section 2.3.3). We also use the term *offset* to describe the start position of a thread; thus the offset of a thread T is str(T).

The first minimal thread $S_1 = T_p$ is at offset s. This means that $T_{2p}$ is at offset 2s; so if 2s < n, $S_2 = T_{2p}$, but if $2s \geq n$, $S_2 = T_{2p+1}$.

## *1.1  s = 1*

We first dispose of the case s = 1:

$str(T_p) = a_2 - pn = s = 1$, and since $len(T_p) = n$, $end(T_p) = n$; but $str(T_0) = 0$, and $end(T_0) = n$: so $T_p$ is covered by $T_0$ and so this case cannot be a stride generator.

## *1.2  s = n/2*

If $s = n/2$, then $T_{2p+1}$ is at offset 0, is of length n-1, and so is completely covered by $T_0$. Since $s \geq 2$, $n \geq 4$, and so $a_2 \geq 4p + 2 = 2(2p + 1)$. So when $s = n/2$, there exists a thread of order $2p+1 \leq a_2/2$ that is covered by $T_0$, which means that $q < a_2/2$. But $n < a_2/p \leq a_2/2$, so $n + q < a_2$ as required.

## 1.3  s < n/2

The threads $S_i$ satisfy $str(S_0) = 0$, $str(S_1) = s$, ... $str(S_k) = ks$, $str(S_{k+1}) = (k+1)s - n$, and so on; so these threads fall into groups, each of which starts with a thread whose offset is a local minimum. We denote these *minimal* threads as $M_0 = S_0 = T_0$, $M_1 = S_{k+1}$, $M_2$ ... Note that $k \geq 2$, since we have assumed that $s < n/2$.

When $t = 0$, $M_1 = S_k$ with offset 0, and so $S_k = T_{kp+1}$ is covered by $T_0$. Now:
    $a_2 = pn + s = pks + pt + s = s(pk + 1) + pt \geq 2(pk + 1) + pt$    since $s \geq 2$

So $kp + 1 \leq a_2/2$, and there exists a thread covered by $T_0$ whose order is $\leq a_2/2$: so $q < a_2/2$, $n < a_2/2$ $\Rightarrow$ $n + q < a_2$ as required. This means that we may henceforth assume that $t \geq 1$.

Each group of threads $S_i$ contains either $k$ or $k + 1$ threads, so the difference in order between consecutive minimal threads is either $kp+1$ or $(k+1)p+1$. The length of each thread $S_i$ is one less than its predecessor $S_{i-1}$, and so the difference in length between consecutive minimal threads is either $k$ or $k+1$. Thus we have established the following bounds for the jth minimal thread $M_j$:
    $str(M_j) < s$
    $ord(M_j) \leq ((k+1)p + 1)j$
    $len(M_j) \leq (n+1) - (k+1) - k(j-1)$      since $len(M_1) = len(T_{kp+1}) = len(T_0) - k - 1 = (n+1) - (k+1)$, and subsequent threads get smaller by at least k each time

For $M_j$ to be covered by $T_0$, we require that $len(M_j) + str(M_j) \leq n+1$, and this will certainly be true if:
    $(n+1) - (k+1) - k(j-1) + (s-1) \leq n+1$ $\Rightarrow$ $k(j-1) \geq s-k-2$ $\Rightarrow$ $j-1 \geq (s-k-2)/k = (s-2)/k - 1$ $\Rightarrow$ $j \geq (s-2)/k$      - (1)
    [Aside: there is no need to show that $len(M_j) \geq 1$, since if $M_j$ does not exist then q is limited in exactly the same way as when $M_j$ exists and is covered by $T_0$].

Clearly there exists an integer $j_0$ satisfying $(s-2)/k + 1 > j_0 \geq (s-2)/k$, and in this case:
    $ord(M_{j_0}) \leq ((k+1)p + 1)((s-2)/k + 1)$      - (2)
We also have:
    $a_2 = pn + s = p(ks + t) + s$ $\Rightarrow$ $a_2 = (kp+1)s + pt$      - (3)

### 1.3.1  k > 2

For the following sections we assume that $k > 2$; in fact, this assumption is required only when $s > k+2$, but the case $k = 2$ is not sensitive to this distinction.

#### 1.3.1.1  p = 2

For this case, we show that $ord(M_{j_0}) \leq a_2/2 + 1$ $\Rightarrow$ $q \leq a_2/2$; since $n < a_2/2$, $n + q < a_2$ follows immediately.
We write $d = a_2/2 + 1 - ord(M_{j_0})$ $\Rightarrow$ $2d = a_2 + 2 - 2*ord(M_{j_0})$; we must show that $d \geq 0$.

From (3) we have:
    $a_2 = (2k+1)s + 2t \geq (2k+1)s + 2$     since $t \geq 1$
From (2) we have:
    $ord(M_{j_0}) \leq (2(k+1) + 1)((s-2+k)/k) = (2k+3)(s+k-2)/k$
So:
    $2d \geq (2k+1)s + 4 - 2(2k+3)(s+k-2)/k$
$\Rightarrow$ $2kd \geq (2k+1)sk + 4k - 2(2k+3)(s+k-2) = s(k(2k+1) - 2(2k+3)) + 4k - 2(2k+3)(k-2) = s(2k^2 - 3k - 6) + 4k - 4k^2 + 2k + 12$
$\Rightarrow$ $2kd \geq s(2k^2 - 3k - 6) - (4k^2 - 6k - 12)$      - (4)

##### 1.3.1.1.1  s ≤ k+2

When $s \leq k+2$, $j_0 = 1$ and we have:
    $ord(M_1) = (k+1)p + 1 = 2k+3$
    $a_2 \geq (2k+1)s + 2$
So:     $2d \geq (2k+1)s + 4 - 4k - 6 = 2(s-2)k + s - 2 \geq 0$     since $s \geq 2$

##### 1.3.1.1.2  s > k+2

From {4}:
    $2kd \geq (k+3)(2k^2 - 3k - 6) - (4k^2 - 6k - 12)$
When $k \geq 3$, $2k^3$ increases faster than $4k^2$ as k increases, and so the first term increases faster than the second. When $k = 3$ we have $2kd \geq 12$, so we have $2kd \geq 0$ for all $k \geq 3$ as required. Note that when $k = 2$ this bound is inadequate: we have only that $2kd \geq -12$.

#### 1.3.1.2  p ≥ 3

Here we show that $ord(M_{j_0}) \leq 2a_2/3$ $\Rightarrow$ $q < 2a_2/3$; since $n < a_2/p \leq a_2/3$, $n + q < a_2$ follows immediately.
We write $d = 2a_2/3 - ord(M_{j_0})$ $\Rightarrow$ $3d = 2a_2 - 3*ord(M_{j_0})$; we must show that $d \geq 0$.

From (3) we have:
    $a_2 \geq (pk+1)s + p$    since $t \geq 1$
From (2) we have:
    $ord(M_{j_0}) \leq ((k+1)p + 1)(s-2+k)/k$
So:
    $3d \geq 2s(pk+1) + 2p - 3((k+1)p + 1)(s+k-2)/k$
$\Rightarrow$ $3kd \geq 2ks(kp+1) + 2kp - 3((k+1)p + 1)(s+k-2) = s(2k(kp+1) - 3((k+1)p+1)) + 2kp - 3((k+1)p+1)(k-2) = As + B$ where:
        $A = p(2k^2 - 3k - 3) + (2k - 3)$
        $B = p(2k - 3(k-2)(k+1)) - 3(k-2) = p(2k - 3k^2 + 3k + 6) - 3(k-2) = -p(3k^2 - 5k - 6) - 3(k-2)$
$\Rightarrow$ $3kd \geq s(p(2k^2 - 3k - 3) + (2k - 3)) - p(3k^2 - 5k - 6) - 3(k - 2)$      - (5)

##### 1.3.1.2.1  s ≤ k+2

When $s \leq k+2$, $j_0 = 1$ and we have:
    $ord(M_1) = (k+1)p + 1$
    $a_2 \geq (pk+1)s + p$

So:   $3d \geq 2(pk+1)s + 2p - 3((k+1)p + 1) = p(2sk + 2 - 3(k+1)) + 2s - 3 \geq p(4k + 2 - 3k - 3) + 1$   since $s \geq 2$
   => $3d \geq p(k-1) + 1 \geq 0$   as required, since $k \geq 2$.

*1.3.1.2.2  s > k+2*
Substituting $s = k+3$ in (5) we have:
   $3kd \geq (k+3)(p(2k^2 - 3k - 3) + (2k - 3)) - p(3k^2 - 5k - 6) - 3(k - 2)$
   => $3kd \geq p((k+3)(2k^2 - 3k - 3) - (3k^2 - 5k - 6)) + (k+3)(2k-3) - 3(k-2) = p(2k^3 - 7k - 3) + (2k^2 - 3) = Xp + y$ where:
       $X > 0$ for all $k \geq 3$, since $2k^3$ increases faster than $7k$ for $k \geq 3$, and $X = 30$ for $k = 3$
       $Y > 0$ for all $k \geq 3$, since $2k^2$ is greater than 3 for $k \geq 3$, and $Y = 15$ for $k = 3$
   => $3kd > 0$   provided $k \geq 3$.  (Note that when $k = 2$, $3kd \geq 5 - p$, which is an insufficient bound when $p > 5$)

*1.3.2  k = 2*
As an example of a stride generator with $k = 2$, $s < n/2$ where the contsraints of 1.x.x above are inadequate, consider $A = \{1, 52, 73\}$ where $a_2 = 52$, $n = 21$, $p = 2$, $s = 10$, $k = 2$ and $t = 1$. (X) shows that we need $j \geq (s-2)/k = 4$ in order to be certain that the thread $M_j$ is covered by $T_0$; but then (Y) guarantees only that $ord(M_j) \leq 28$, whereas our argument requires $ord(M_j) \leq a_2/2 + 1 = 27$. In practice, of course, all is well; both $M_3$ and $M_4$ are covered by $T_0$:

| i | $ord(S_i)$ | $str(S_i)$ | $len(S_i)$ | $str(S_i)+len(S_i)$ | |
|---|---|---|---|---|---|
| 1 | 2 | 10 | 21 | 31 | |
| 2 | 4 | 20 | 20 | 40 | |
| 3 | 7 | 9 | 19 | 28 | $j = 1$ |
| 4 | 9 | 19 | 18 | 37 | |
| 5 | 12 | 8 | 17 | 25 | $j = 2$ |
| 6 | 14 | 18 | 16 | 34 | |
| 7 | 17 | 7 | 15 | 22 | $j = 3$: (just) covered by $T_0$ |
| 8 | 19 | 17 | 14 | 31 | |
| 9 | 22 | 6 | 13 | 19 | $j = 4$: covered by $T_0$ |

We have seen above that we need to consider this case separately only when $s > k+2$ => $s \geq 5$.

Each group of threads $\{S_i\}$ contains either k or (k+1) threads, and the difference in order and length between consecutive minimal threads is determined accordingly:

|   | $ord(M_j) - ord(M_{j-1})$ | $len(M_j) - len(M_{j-1})$ |
|---|---|---|
| k | kp + 1 | k |
| k + 1 | (k + 1)p + 1 | k + 1 |

So precise formulae are:
   $len(M_j) = (n+1) - j_1(k+1) - j_2 k$          for some  $j_1 \geq 1$, $j_2 \geq 0$  with $j_1 + j_2 = j$          - (6)
   $ord(M_j) = ((k+1)p + 1)j_1 + (kp+1)j_2$          - (7)
Substituting $k = 2$, we have:
   $len(M_j) = (n+1) - (3j_1 + 2j_2)$          - (8)
   $ord(M_j) = (3p + 1)j_1 + (2p + 1)j_2 = p(3j_1 + 2j_2) + (j_1 + j_2) \leq p(3j_1 + 2j_2) + (3j_1 + 2j_2)/2 = (2p+1)(3j_1 + 2j_2)/2$          - (9)
$M_{j0}$ is certainly covered by $T_0$ when $len(M_{j0}) + (s-1) \leq n+1$ => $s-1 \leq 3j_1 + 2j_2$. Since $j_1, j_2$ are both integers, we know that if $j_1$ and $j_2$ are the smallest values which satisfy $3j_1 + 2j_2 \geq s-1$, then $3j_1 + 2j_2 < s+2$; so from (9) we have:
   $ord(M_{j0}) \leq (2p+1)(s+1)/2$          - (10)

*1.3.2.1  p = 2*
We know that $n < a_2/2$, and so we have only to show that $ord(M_{j0}) \leq a_2/2 + 1$ => $q \leq a_2/2$ => $n + q < a_2$.
From (10) we have:
   $ord(M_{j0}) \leq (5s + 5)/2$          - (11)
and we have:
   $a_2 = pn + s = (kp + 1)s + pt = 5s + 2t \geq 5s + 2$          - (12)

*1.3.2.1.1  s even*
We write $s = 2u$; from (11) we get:
   $ord(M_{j0}) \leq (10u + 5)/2 = 5u + (5/2)$ => $ord(M_{j0}) \leq 5u + 2$   since $ord(M_{j0})$ is integral
From (12) we have:
   $a_2/2 \geq (10u + 2)/2 = 5u + 1$ => $a_2/2 + 1 \geq 5u + 2$
So  $ord(M_{j0}) \leq a_2/2 + 1$  as required.

*1.3.2.1.2  s odd*
We write $s = 2u + 1$  ($u \geq 0$); from (11) we have:
   $ord(M_{j0}) \leq 5u + 5$
and from (12):
   $a_2/2 \geq 5u + (7/2) > 5u + 3$ => $a_2/2 + 1 > 5u + 4$
This leaves the possibility that when $t = 1$ $ord(M_{j0}) = a_2/2 + 2$; but for $t > 1$ the result is proved. There are two ways to deal with $t = 1$:
   a) It is easy to see that when $t = 1$, we have $j_1 = 1$, and $j_2 = (j-1)$:

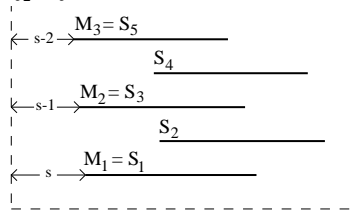

      Substituting in our original exact formulae (6) and (7) we have:
         $ord(M_j) = 7 + 5(j_0 - 1) = 5j + 2$
         $len(M_j) = (n+1) - 3 - 2(j-1) = (n+1) - 2j - 1$

So $M_j$ is covered by $T_0$ when $s-1 \leq 2j+1 \Rightarrow j \geq (s-2)/2$; and we can certainly find such a $j$ where $j < s/2$.
When $j < s/2$, $\text{ord}(M_j) < 5s/2 + 2 = (10u+5)/2 + 2 = 5u + 9/2$; but $\text{ord}(M_j)$ is integral, so $\text{ord}(M_j) \leq 5u + 4$ as required.

b) A simpler approach is to consider the sum $n + q$ as a whole for the particular case $s = 2u + 1$, $t = 1$ (but this is a weaker result than (a), since it does not prove that $q \leq a_2/2$):
We have:
$$n = 2s + 1 = 4u + 3$$
$$a_2 = 2n + s = 8u + 6 + 2u + 1 = 10u + 7$$
$$q \leq \text{ord}(M_{j0}) - 1 = 5u + 4$$
So: $n + q \leq 10u + 7 = a_2$ as required.

### 1.3.2.2 $p \geq 3$

We know that $n < a_2/3$, and will show that $\text{ord}(M_{j0}) < 2a_2/3 \Rightarrow n + q < a_2$.
(10) gives:
$$3*\text{ord}(M_{j0}) \leq 3(2p+1)(s+1)/2$$
and:
$$2a_2 = 2(pn + s) = 2((kp+1)s + pt) \geq 2((2p+1)s + p) \quad \text{since } t \geq 1$$
So:
$$2a_2 - 3*\text{ord}(M_{j0}) \geq 2((2p+1)s + p) - 3((2p+1)s + (2p+1))/2 = s(2p+1)/2 + 2p - 3p - (3/2) \geq 5(2p+1)/2 - p - (3/2) \text{ since } s \geq 5$$
$$= 4p + 1 > 0$$
So $2a_2 > 3*\text{ord}(M_{j0}) \Rightarrow \text{ord}(M_{j0}) < 2a_2/3$ as required.

## 1.4 $s > n/2$

In this case, the threads over $T_0$ group as follows (cf diagram in section 2.3.3 above, which illustrates the case for $s < n/2$):

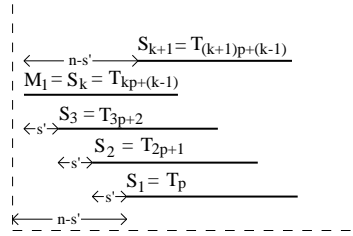

Each thread $S_i$ is one shorter than its predecessor, and the difference in order is either $p+1$ (when going 'up' the staircase from the right) or $p$ (when following a minimal thread $M_i$). We write $s' = n-s \Rightarrow s = n-s'$, $0 < s' < n/2$; we also treat $s' = 1$ as a special case, and so can assume $s' \geq 2$ for the general case. We have:
$$a_2 = pn + s = (p+1)n - s'$$
and we write:
$$n = ks' + t \quad 0 \leq t < s', \; k \geq 2$$
If $t = 0$, $M_1$'s offset is zero, and so $M_1$ is completely covered by $T_0$; we also have:
$$\text{ord}(M_1) = kp + k - 1$$
and: $a_2 = (p+1)n - s' = (p+1)ks' - s' = (kp + k - 1)s'$
Since $s' \geq 2$, this means that $\text{ord}(M_1) \leq a_2/2$, and so $q < a_2/2$; since $n < a_2/2$, this gives $n + q < a_2$ as required for $t = 0$; so we may assume that $t \geq 1$ in what follows. [A limiting example where $t = 0$ and $\text{ord}(M_1) = a_2/2$ is given by $\{1, 30, 38\}$ where $s' = 2$, $n = 8$, $p = 3$, $k = 4$, and $M_1 = T_{15}$].

### 1.4.1 $s' = 1$

When $s' = 1$, $M_1 = S_n = T_{np+(n-1)}$, so $\text{ord}(M_1) = a_2$ (since $a_2 = pn + s = pn + (n-1)$); this gives us a bound for $q$ of $q < a_2$, which is certainly not good enough to show that $n + q < a_2$. Instead, we must look at the threads $S_i$ themselves:
$$\text{ord}(S_i) = (p+1)i - 1$$
$$\text{str}(S_i) = n - i$$
$$\text{len}(S_i) = (n+1) - i$$
So $S_i$ is covered by $T_0$ as soon as $\text{str}(S_i) + \text{len}(S_i) \leq n+1 \iff n - i + (n+1) - i \leq (n+1) \iff 2i \geq n$; we consider three cases: $n$ even, $n = 1$, $n$ odd.

#### 1.4.1.1 $n$ even

We write $n = 2m$, and consider $S_m$ which is (just) covered by $T_0$.
$$\text{ord}(S_m) = (p+1)m - 1$$
$$a_2 = (p+1)n - s' = 2(p+1)m - 1$$
So $\text{ord}(S_m) < a_2/2 \Rightarrow q < a_2/2 \Rightarrow n + q < a_2$ as required.

#### 1.4.1.2 $n$ odd, $n \geq 3$

We write $n = 2m-1$, and consider $S_m$ which is - again - just covered by $T_0$.
$$\text{ord}(S_m) = (p+1)m - 1$$
$$a_2 = (p+1)n - s' = (2m-1)(p+1) - 1$$

##### 1.4.1.2.1 $p = 2$

We show $\text{ord}(S_m) \leq a_2/2 + 1 \Rightarrow q \leq a_2/2$; since $p = 2$, $n < a_2/2$, and so $n + q < a_2$ as required.
$$\text{ord}(S_m) = 3m - 1$$
$$a_2 = 3(2m-1) - 1 = 6m - 4$$
So $\text{ord}(S_m) = a_2/2 + 1$, which is just sufficient.

##### 1.4.1.2.2 $p \geq 3$

We show $\text{ord}(S_m) < 2a_2/3 \Rightarrow q < 2a_2/3$; since $p \geq 3$, $n < a_2/3$, and so $n + q < a_2$ as required.
$$2a_2 - 3*\text{ord}(S_m) = 2(2m-1)(p+1) - 2 - 3(p+1)m + 3 = (p+1)(m-2) + 1 > 0 \text{ for all } m \geq 2.$$
So $\text{ord}(S_m) < 2a_2/3$ (and the result is proved) provided that $m > 1$; this leaves just $m = 1 \iff n = 1$ to consider.

#### 1.4.1.3 $n = 1$

When $n = 1$, we have $A = \{1, a_2, a_2+1\}$ and the thread diagram for the whole of the range $0 \leq x < a_2$ looks like this:

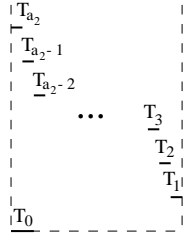

Each thread $T_i$ for i = 1 ... $a_2$ is of length 1, and $str(T_i) = (i+1)a_3 - ia_3 = a_2 - i$. $T_0$ is of length 2, and so $T_{a(2)-1}$ is the first thread which is covered by $T_0$.
So $q < a_2-1$, and n = 1 => n + q < $a_2$. [In fact, we can easily see that such stride generators are canonical: the diagram above shows that p = $a_2$-2, and so we have shown that $T_{p+1}$ is completely covered by $T_0$, and so the stride generator is canonical].
[This case is really *not* a sub-case of s' = 1: n = 1 => $a_2$ = pn + 0, which turns out not to be a stride generator at all. Maybe some more thought should go into this to get it really straight!]

*1.4.2  s' ≥ 2*

Referring to the diagram in 1.3 above, each group of threads $\{S_i\}$ has either k or k+1 threads in it, with the first group containing just k threads. We therefore obtain the following bounds on minimal threads $M_j$:

$str(M_j) < s'$
$ord(M_j) \leq k(p+1) - 1 + ((k+1)(p+1) - 1)(j-1)$
$len(M_j) \leq (n+1) - kj$

For $M_j$ to be covered by $T_0$ we need $len(M_j) + str(M_j) \leq n+1$, and this will certainly be true when:

$(n+1) - kj + s' - 1 \leq n+1 \iff kj \geq s'-1 \iff j \geq (s'-1)/k$   - (13)

Clearly there is an integer $j_0$ satisfying $(s'-1)/k + 1 > j_0 \geq (s'-1)/k$, and then we have:

$ord(M_{j_0}) \leq k(p+1) - 1 + ((k+1)(p+1) - 1)(s'-1)/k$   - (14)

We also have:

$a_2 = (p+1)n - s' = (p+1)(ks' + t) - s' = s'(k(p+1) - 1) + (p+1)t$   - (15)

*1.4.2.1  k > 2*

We see below that it is only for the case s' > k+1 that we need to assume k ≥ 3.

*1.4.2.1.1  s' ≤ k+1*

When s' < k+1, (13) shows that $j_0$ = 1 is sufficient to ensure that $M_{j_0} = M_1$ is covered by $T_0$. We now show that $ord(M_1) < a_2/2 \Rightarrow q < a_2/2$; since $n < a_2/2$, this shows $n + q < a_2$ as required.

$ord(M_1) = k(p+1) - 1$
$a_2 = s'(k(p+1) - 1) + (p+1)t \geq s'(k(p+1) - 1) + (p+1)$ since t ≥ 1
$> s' * ord(M_1) \geq 2 * ord(M_1)$  since s' ≥ 2

So $ord(M_1) < a_2/2$ as required.

*1.4.2.1.2  s' > k+1*

Again, we consider the cases p = 2 and p ≥ 3 separately.

*1.4.2.1.2.1  p = 2*

We show $ord(M_{j_0}) \leq a_2/2 + 1 \Rightarrow q \leq a_2/2$; since $n < a_2/2$, this shows $n + q < a_2$ as usual.
We write $d = a_2/2 + 1 - ord(M_{j_0})$; from (14) and (15) we have:

$2d \geq s'(k(p+1) - 1) + (p+1)t + 2 - 2k(p+1) + 2 - 2((k+1)(p+1) - 1)(s'-1)/k$
$\Rightarrow 2kd \geq ks'(k(p+1) - 1) + k(p+1) + 4k - 2k^2(p+1) - 2((k+1)(p+1) - 1)(s'-1)$   since t ≥ 1
$= ks'(3k-1) + 3k + 4k - 6k^2 - 2((k+1)3 - 1)(s'-1)$   since p = 2
$= s'(3k^2 - 7k - 4) - (6k^2 - 13k - 4)$

Now $k \geq 3 \Rightarrow s' \geq 5$, and $3k^2 - 7k - 4 > 0$; so:
$2kd \geq 5(3k^2 - 7k - 4) - (6k^2 - 13k - 4) = 9k^2 - 22k - 16 \geq -1$ for all $k \geq 3$.

So $d \geq -(1/2k) \geq -1/6$ for all $k \geq 3$, and so $ord(M_{j_0}) \leq a_2/2 + 1 + (1/6)$; but $ord(M_{j_0})$ is integral, and so $ord(M_{j_0}) \leq a_2/2 + 1$ as required.
[ This is a 'sharp' bound: consider k = 3, s' = 5, t = 1, p = 2 => n = 16, $a_2$ = 43, $a_3$ = 59; we find $j_0$ = 2 which gives $ord(M_{j_0}) \leq 22+(2/3)$ and hence $ord(M_{j_0}) \leq 22$ as necessary.]

*1.4.2.1.2.2  p ≥ 3*

We show $ord(M_{j_0}) < 2a_2/3 \Rightarrow q < 2a_2/3$; since $p \geq 3$, $n < a_2/3$ and we have $n + q < a_2$.
We write $d = 2a_2/3 - ord(M_{j_0})$; from (14) and (15) we have:

$3d \geq 2s'(k(p+1) - 1) + 2(p+1)t - 3k(p+1) + 3 - 3((k+1)(p+1) - 1)(s'-1)/k$
$\Rightarrow 3kd \geq 2ks'(k(p+1) - 1) + 2k(p+1) - 3k^2(p+1) + 3k - 3((k+1)(p+1) - 1)(s'-1)$
$= s'(2k(k(p+1) - 1) - 3((k+1)(p+1) - 1)) + (2k(p+1) - 3k^2(p+1) + 3k + 3((k+1)(p+1) - 1))$
$= s'(p(2k^2 - 3k - 3) + 2k^2 - 2k - 3k) + (p(2k - 3k^2 + 3k + 3) + 2k - 3k^2 + 3k + 3k)$
$= s'(p(2k^2 - 3k - 3) + k(2k - 5)) - (p(3k^2 - 5k - 3) + k(3k - 8))$

Now $s' \geq k+2$, and $p(2k^2 - 3k - 3) + k(2k - 5) > 0$ since $p > 0$ and $k \geq 3$, so:

$3kd \geq (k+2)(p(2k^2 - 3k - 3) + k(2k - 5)) - (p(3k^2 - 5k - 3) + k(3k - 8))$
$= p((k+2)(2k^2 - 3k - 3) - 3k^2 + 5k + 3) + (k(k+2)(2k - 5) - k(3k - 8))$
$= p(2k^3 - 2k^2 - 4k - 3) + k(2k^2 - 4k - 2) > 0$  for all $k \geq 3$

So $3kd > 0 \Rightarrow d > 0 \Rightarrow ord(M_{j_0}) < 2a_2/3$ as required.

*1.4.2.2  k = 2*

When s' > k+1, the bounds above are not sufficient to demonstrate $n + q < a_2$ when k = 2: for p = 2 we obtain only $2kd \geq 6(1-s')$, and for $p \geq 3$ we have only $3kd \geq (p+2)(1-s')$. Instead we must develop improved bounds by working with precise formulae for the length and order of the threads $M_j$.
Each group of threads $\{S_i\}$ contains either k or (k+1) threads, and the difference in order and length between consecutive minimal threads is determined accordingly:

| | $ord(M_j) - ord(M_{j-1})$ | $len(M_j) - len(M_{j-1})$ |
|---|---|---|
| k | $k(p+1) - 1$ | k |
| k + 1 | $(k+1)(p+1) - 1$ | k + 1 |

So precise formulae are:

$\quad\quad\quad$ len($M_j$) = (n+1) - $kj_1$ - (k+1)$j_2$ $\quad\quad\quad\quad\quad\quad$ for some $j_1 \geq 1$, $j_2 \geq 0$ with $j_1 + j_2 = j$ $\quad\quad$ - (16)
$\quad\quad\quad$ ord($M_j$) = (k(p+1) - 1)$j_1$ + ((k+1)(p+1) - 1)$j_2$ $\quad\quad\quad\quad\quad\quad\quad\quad\quad\quad\quad\quad\quad\quad\quad\quad\quad\quad$ - (17)
For k = 2 we have:
$\quad\quad\quad$ len($M_j$) = (n+1) - (2$j_1$ + 3$j_2$) $\quad\quad\quad\quad\quad\quad\quad\quad\quad\quad\quad\quad\quad\quad\quad\quad\quad\quad\quad\quad\quad\quad\quad\quad$ - (18)
and $\quad\quad$ ord($M_j$) = (2p+1)$j_1$ + (3p+2)$j_2$ $\quad\quad\quad\quad\quad\quad\quad\quad\quad\quad\quad\quad\quad\quad\quad\quad\quad\quad\quad\quad\quad\quad$ - (19)
$\quad\quad\quad\quad\quad\quad$ = p(2$j_1$ + 3$j_2$) + ($j_1$ + 2$j_2$) $\leq$ p(2$j_1$ + 3$j_2$) + (2/3)(2$j_1$ + 3$j_2$) = (p + 2/3)(2$j_1$ + 3$j_2$)
$M_{j0}$ is certainly covered by $T_0$ when len($M_{j0}$) + (s'-1) $\leq$ (n+1) $\iff$ s'-1 $\leq$ (2$j_1$ + 3$j_2$); since $j_1$, $j_2$ are integral, we know that s'+2 > (2$j_1$ + 3$j_2$) $\geq$ s'-1, and so 2$j_1$ + 3$j_2$ $\leq$ s'+1; so:
$\quad\quad\quad$ ord($M_{j0}$) $\leq$ (p + 2/3)(s'+1) $\quad\quad\quad\quad\quad\quad\quad\quad\quad\quad\quad\quad\quad\quad\quad\quad\quad\quad\quad\quad\quad\quad\quad\quad\quad\quad$ - (20)

*1.4.2.2.1  p = 2*
From (20) and (15) we have:
$\quad\quad\quad$ ord($M_{j0}$) $\leq$ (2 + 2/3)s' + (2 + 2/3)
$\quad\quad\quad$ $a_2$ = s'(k(p+1)-1) + (p+1)t $\geq$ s'(2p+1) + (p+1) = 5s' + 3 $\quad\quad$ since t $\geq$ 1
These bounds are not sufficient to show that ord($M_{j0}$) $\leq$ $a_2$/2 + 1  for all s', so we must split into separate cases again.

*1.4.2.2.1.1  t $\geq$ s'/3*
We have:
$\quad\quad\quad$ $a_2$ $\geq$ 5s' + 3(s'/3) = 6s'
So:
$\quad\quad\quad$ $a_2$/2 + 1 - ord($M_{j0}$) $\geq$ s'/3 + 1 - (2 + 2/3) = s'/3 - 5/3 $\geq$ -1/3 $\quad\quad$ since s' > k+1 => s' $\geq$ 4
So ord($M_{j0}$) $\leq$ $a_2$/2 + 1 + 1/3, which, since ord($M_{j0}$) is integral, is sufficient to show that ord($M_{j0}$) $\leq$ $a_2$/2 + 1 => q $\leq$ $a_2$/2; since n < $a_2$/2, n + q < $a_2$ follows immediately.

*1.4.2.2.1.2  t < s'/3*
The next section contains an argument to show that q $\leq$ $a_2$/2 in this case, too, but a simpler approach is to consider the sum n + q directly:
$\quad\quad\quad$ n = ks' + t = 2s' + t < (2 + 1/3)s'
$\quad\quad\quad$ q < ord($M_{j0}$) $\leq$ (2 + 2/3)s' + (2 + 2/3)
$\quad$ => n + q < 5s' + (2 + 2/3) < 5s' + 3 = $a_2$ as required.

*1.4.2.2.1.3  Alternative method for all t*
From (18), (19) and (15) we have the following exact formulae:
$\quad\quad\quad$ len($M_j$) = (n+1) - (2$j_1$ + 3$j_2$) = (n+1) - (2j + $j_2$) $\quad\quad\quad\quad\quad\quad\quad\quad\quad\quad\quad\quad\quad\quad\quad$ - (21)
$\quad\quad\quad$ ord($M_j$) = 5$j_1$ + 8$j_2$ = 5j + 3$j_2$ $\quad\quad\quad\quad\quad\quad\quad\quad\quad\quad\quad\quad\quad\quad\quad\quad\quad\quad\quad\quad\quad\quad\quad$ - (22)
$\quad\quad\quad$ $a_2$ = 5s' + 3t $\quad\quad\quad\quad\quad\quad\quad\quad\quad\quad\quad\quad\quad\quad\quad\quad\quad\quad\quad\quad\quad\quad\quad\quad\quad\quad\quad\quad\quad\quad$ - (23)
where $j_2$ is the number of groups of threads which have 3 - rather than 2 - members.
Now from (21) $M_{j0}$ is certainly covered by $T_0$ when  len($M_{j0}$) + (s'-1) $\leq$ (n+1) $\iff$ 2j + 3$j_2$ $\geq$ (s'-1) $\iff$ 2j $\geq$ (s'-1) - $j_2$. This is certainly true if:
$\quad\quad\quad$ j $\geq$ (s'-1)/2 $\quad\quad\quad\quad\quad\quad\quad\quad\quad\quad\quad\quad\quad\quad\quad\quad\quad\quad\quad\quad\quad\quad\quad\quad\quad\quad\quad\quad\quad\quad\quad$ - (24)
The thread groups are themselves grouped as follows:

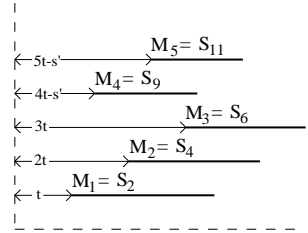

Since n = 2s' + t => n - 2s' = t, we have:
$\quad\quad\quad$ str($S_1$) = n-s'
$\quad\quad\quad$ str($S_2$) = n-2s' = t $\quad\quad\quad\quad\quad\quad\quad\quad\quad\quad\quad\quad$ this is $M_1$ (a group of 2)
$\quad\quad\quad$ str($S_3$) = (n-3s') + n = 2n-3s'
$\quad\quad\quad$ str($S_4$) = 2n-4s' = 2t $\quad\quad\quad\quad\quad\quad\quad\quad\quad\quad\quad\quad$ this is $M_2$ (a group of 2)
$\quad\quad\quad$ ....
$\quad\quad\quad$ str($S_{2v}$) = vn - (2v)s' = vt > s' $\quad\quad$ where s' = vt - w, 0 $\leq$ w < t
$\quad\quad\quad$ str($S_{2v+1}$) = vn - (2v+1)s' = vt-s' $\quad\quad\quad\quad\quad\quad\quad\quad$ this is $M_v$ (the first group of 3)
$\quad\quad\quad$ ....
We see that:
$\quad\quad\quad$ str($M_1$) = t
$\quad\quad\quad$ str($M_2$) = 2t
$\quad\quad\quad\quad$ ....
$\quad\quad\quad$ str($M_j$) = jt mod s'
and each time s' has to be subtracted, a 3-thread group is present.
So: $\quad\quad$ If jt = fs' + g,  0 $\leq$ g < s', then the minimal threads $M_1$, $M_2$, ... $M_j$ include exactly f 3-thread groups. $\quad\quad$ - (25)

We now consider two sub-cases - s' even, and s' odd - separately, and show ord($M_{j0}$) $\leq$ $a_2$/2 => q < $a_2$/2; but n < $a_2$/2 => n + q < $a_2$ as required.
$\quad\quad$ a) s' even, s' = 2u:
$\quad\quad\quad\quad$ From (24), $j_0$ = u will ensure that $M_{j0}$ is covered by $T_0$.
$\quad\quad\quad\quad$ From (25), there are exactly f 3-thread groups in $M_1$ ... $M_{j0}$, where ut = 2uf + g, 0 $\leq$ g < 2u
$\quad\quad\quad\quad\quad$ => t = 2f + g/u => f = t/2 - g/2u => f $\leq$ t/2
$\quad\quad\quad\quad$ But f is precisely the $j_2$ of (22), so we have:
$\quad\quad\quad\quad\quad\quad$ ord($M_{j0}$) $\leq$ 5$j_0$ + 3t/2 = 5u + 3t/2 => 2*ord($M_{j0}$) $\leq$ 10u + 3t = 5s' + 3t = $a_2$ (see (23))
$\quad\quad\quad\quad$ So ord($M_{j0}$) $\leq$ $a_2$/2.
$\quad\quad$ b) s' odd, s' = 2u+1:
$\quad\quad\quad\quad$ From (24), $j_0$ = u will ensure that $M_{j0}$ is covered by $T_0$.
$\quad\quad\quad\quad$ From (25), there are exactly f 3-thread groups in $M_1$ ... $M_{j0}$, where ut = f(2u+1) + g,  0 $\leq$ g < 2u+1

$\Rightarrow \quad f = (ut - g)/(2u+1) < (ut - g)/2u = (t/2) - (g/2u) < t/2$

But f is precisely the $j_2$ of (22), so we have:

$\text{ord}(M_{j0}) \leq 5j_0 + 3t/2 \Rightarrow 2*\text{ord}(M_{j0}) \leq 10j + 3t = 10u + 3t$, and $a_2 = 5s' + 3t = 10u + 3t + 5$

So $\text{ord}(M_{j0}) \leq a_2/2$.

There are two interesting points to note about this argument.

a) The argument does not rely on the exact formula for $\text{len}(M_j)$ given in (21); that is, $j_2$ is discounted. $j_2$ is a count of the number of "3-thread" groups present in $M_1 \ldots M_{j0}$, and it is not clear whether there is always at least one such group. However, here is an example where the bound used for $\text{len}(M_{j0})$ is 'sharp*': {1, 56, 78}, where n = 22, s' = 10, t = 2:

| j | $\text{ord}(M_j)$ | $\text{len}(M_j)$ | (21) |
|---|---|---|---|
| 0 | 0 | 23 | 23 |
| 1 | 5 | 21 | 21 |
| 2 | 10 | 19 | 19 |
| 3 | 15 | 17 | 17 |
| 4 | 20 | 15 | 15 |
| $j_0 = 5$ | 28 | 12 | 13 |

[ * actually, this doesn't seem to be the case ]

b) The argument does not depend on $t < s'/3$ (or $s' \geq 4$), and so we could use it for all $s > n/2$, $k = 2$, $p = 2$.

There is, presumably, a corresponding argument for $s < n/2$.

### 1.4.2.2.2 $p \geq 3$

From (20):

$3*\text{ord}(M_{j0}) \leq (3p+2)(s'+1)$

From (15):

$2a_2 \geq 2s'(2p+1) + 2(p+1)$

So:

$2a_2 - 3*\text{ord}(M_{j0}) \geq 4s'p + 2s' + 2p + 2 - 3s'p - 3p - 2s' - 2 = s'p - p = p(s'-1) > 0$ since $s' \geq 4$

So $\text{ord}(M_{j0}) < 2a_2/3 \Rightarrow q < 2a_2/3$; since $p \geq 3$, $n < a_2/3$ and so $n + q < a_2$ as required.

## 2 $C_1 > a_2/2$

Following section 2.4.3 we write $n' = a_2 - C_1$, and note that $pn' < a_2 < (p+1)n'$; we write:

$a_3 = a_2 + C_1 = a_2 + n = 2a_2 - n'$
$a_2 = pn' + s \quad 1 \leq s < n$
$n' = ks + t \quad 0 \leq t < s$
$s' = n' - s$

Experiment suggests that we should be able to show that $q \leq a_2/2$ and $n < a_2/2$, and thus $n + q < a_2/2$. We will see below that we manage $q \leq a_2/2$ for $p = 2$, but not quite for $p \geq 3$; nonetheless, we start by proving $n < a_2/2$ regardless of the value of p.

From section 2.4.3 above we have $n = n' + p - 2$, $pn' < a_2$, so $n < n' + a_2/n' - 2$; we now consider:

$f(x) = x + a_2/x - 2$ over the range $0 < x < a_2/2$

$f'(x) = 0$ when $1 - a_2/x^2 = 0 \Leftrightarrow x = +/- \sqrt{a_2}$; we find $f(x) \to \infty$ as $x \to 0$; $f(a_2/2) = a_2/2 + 2 - 2 = a_2/2$; and $f(\sqrt{a_2}) = 2\sqrt{a_2} - 2$. Now $2\sqrt{a_2} - 2 < a_2/2 \Leftrightarrow 4\sqrt{a_2} < a_2 + 4$ which is true for $a_2 \geq 5$ - and in this case we also have $\sqrt{a_2} < a_2/2$. So for $a_2 \geq 5*$ we have the following shape for the curve f(x):

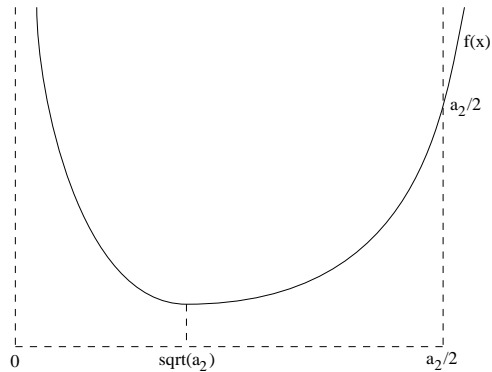

Furthermore, $f(2) = 2 + a_2/2 - 2 = a_2/2$, so:

$f(x) \leq a_2/2$ for $2 \leq x \leq a_2/2$

Now $n' = 1 \Rightarrow a_2 - C_1 = 1 \Rightarrow a_2 - 1 = C_1 \Rightarrow a_2 - C_2 \leq C_1$, and so Lemma 15 tells us that $n' = 1 \Rightarrow$ the stride generator is of order 0, and so is canonical; so we may assume $n' \geq 2$. Since $n < f(n')$, we have proved that $n < a_2/2$ for all $2 \leq n' < a_2/2$ as required; we see that n approaches $a_2/2$ as n' approaches $a_2/2$ and as it approaches 2.

[ * When $a_2 = 4$, the only value $C_1 > a_2/2$ is $C_1 = 3 \Rightarrow n' = 1$; for $a_2 = 3$, $C_1 > a_2/2 \Rightarrow C_1 = 2 \Rightarrow n' = 1$; for $a_2 = 2$ there is no such value $C_1$. ]

### 2.1 $p = 2$

For this case, we prove $q \leq a_2/2$ by demonstrating the existence of a thread $S_i$ covered by $T_0$ with $\text{ord}(S_i) \leq a_2/2 + 1$; since $n < a_2/2$, this proves that $n + q < a_2$ as required.

We have $a_2 = 2n' + s$, $0 < s < n$, and find (see diagram in section 2.4.3) that:

$T_1$ is at offset n', with length n'
$S_1 = T_3$ is at offset n'-s, with length n'-1
If $n'-2s \geq 0$, $S_2 = T_5$ is at offset n'-2s, with length n'-2
If $n'-2s < 0$, $S_2 = T_6$ is at offest 2n'-2s, with length n'-3

This is why we consider the cases $s < n'/2$, $s = n'/2$ and $s > n'/2$ separately.

### 2.1.1 $s = n'/2$

In this case, $S_2 = T_5$ is at offset 0, and so is covered by $T_0$; so $q \leq 4$.
When $s \geq 2$, $n' \geq 4 \implies a_2 = 2n' + s \geq 10$; so for $s \geq 2$, $q < a_2/2$ as required.
When $s = 1$, $n' = 2$ and we find $S_1 = T_3$ at offset 1, length 1 and so $S_1$ is covered by $T_0$ (which is at offset 0 with length 2); so we have $s = 1$, $n' = 2$, $a_2 = 5$, and $q \leq 2$; so $q < a_2/2$ as required.

### 2.1.2  $s < n'/2$

In this case, the threads with offset $< n'$ group as follows (we use $k = 3$ for the example):

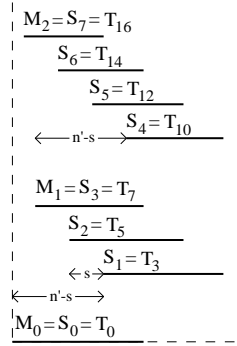

Note that we do not include $T_1$ (since this has offset $= n'$), but it is convenient to think of $M_0 = S_0 = T_0$ as the minimal thread that precedes the first group $\{T_3, T_5, T_7\}$.

The number of threads in each group is determined by the offset of the first thread in the group, which satisfies $n'-s \leq$ offset $< n'$; since $n' = ks + t$, there are either $k$ or $k+1$ threads in each group.

The difference in order between the threads in each group is 2, and the difference in order between the previous minimal thread and the first thread of the next group above it is 3; so the difference in order between successive minimal threads is either $2k+1$ or $2k+3$.

The difference in length between successive threads in the same group is 1, and between the previous minimal thread and the first thread of the next group is 2; so the difference in length between successive minimal threads is either $k+1$ or $k+2$.

The first group always contains $k$ threads (since $str(S_1) = n'-s$), so we have the following bounds:

$\quad ord(M_j) \leq 2k+1 + (j-1)(2k+3)$ \hfill - (26)
$\quad len(M_j) \leq (n'+1) - (k+1) - (j-1)(k+1) = (n'-k) - (j-1)(k+1)$ \hfill - (27)
$\quad str(M_j) < s$ \hfill - (28)

We now show that we may assume $t > 0$; for if $t = 0$, $str(M_1) = 0$ and so $M_1$ is covered by $T_0$; we have:

$\quad ord(M_1) = 2k+1$, and so $q \leq 2k$
$\quad a_2 = 2n' + s = 2ks + 2t + s \geq (2k+1)s$

So $q \leq a_2/2$ as required, provided that $s \geq 2$; we deal with the case $s = 1$ specially below, so we may otherwise assume that $t > 0$.

For $s \geq 2$, we show $q \leq a_2/2$ by finding a minimal thread $M_{j0}$ that is covered by $T_0$ and which satisfies $ord(M_{j0}) \leq a_2/2 + 1$. From (27) and (28), $M_j$ is certainly covered by $T_0$ when:

$\quad len(M_j) + s \leq n'+1 \iff (n'-k) - (j-1)(k+1) + s \leq n'+1 \iff s \leq j(k+1) \iff j \geq s/(k+1)$ \hfill - (29)

We write $d = a_2/2 + 1 - ord(M_j)$:

$\quad d \geq n' + s/2 + 1 - (2k+1) - (j-1)(2k+3) = ks + s/2 + t + 1 - (2k+1) - (j-1)(2k+3) = (2k+1)(s/2 - 1) + t + 1 - (j-1)(2k+3)$ \hfill - (30)

(29) is satisfied by the integer $j_0$ satisfying $s/(k+1) + 1 > j_0 \geq s/(k+1)$; substituting in (30) we have:

$\quad d \geq (2k+1)(s/2 - 1) + t + 1 - (s/(k+1))(2k+3) \geq (2k+1)(s/2 - 1) + 2 - s(2k+3)/(k+1) \quad$ since $t \geq 1$
$\quad\quad = ( (2k+1)(k+1)(s/2 - 1) + 2(k+1) - s(2k+3) )/(k+1)$

So:
$\quad 2(k+1)d \geq (2k^2 + 3k + 1)(s-2) + 4(k+1) - 2s(2k+3) = s(2k^2 + 3k + 1 - 4k - 6) - (4k^2 + 6k + 2 - 4k - 4)$
$\quad\quad\quad\quad\quad = s(2k^2 - k - 5) - (4k^2 + 2k - 2)$ \hfill - (31)

#### 2.1.2.1  $s \geq 18$

Since $s < n'/2$, we have $k \geq 2$, and hence $(2k^2 - k - 5) \geq 0$; so when $s \geq 18$, (31) gives:

$\quad 2(k+1)d \geq 36k^2 - 18k - 90 - 4k^2 - 2k + 2 = 32k^2 - 20k - 88 \geq 0$ for all $k \geq 2$

So for $s \geq 18$, we have $d \geq 0 \implies ord(M_{j0}) \leq a_2/2 + 1 \implies q \leq a_2/2$ as required.
Sadly, for smaller $s$ we must use (29) and (30) directly.

#### 2.1.2.2  $2 \leq s \leq 17$

From (29) we know we must choose $j_0 \geq s/(k+1)$, and so choosing $j_0 \geq s/3$ will be sufficient for all values of $k \geq 2$; we therefore choose $j_0$ to be the smallest integer $\geq s/3$.

From (30) we have:

$\quad d \geq (2k+1)(s/2 - 1) + 2 - (j-1)(2k+3) = (2k+1)(s/2 - 1 - j + 1) + 2 - 2(j-1) = (2k+1)(s/2 - j) - (2j-4)$

From the table below we see that $(s/2 - j_0) \geq 0$ for the values we are considering; since $k \geq 2$ we therefore have:

$\quad d \geq 5(s/2 - j) - (2j-4) = (5s/2) - 7j + 4$

The table shows that $d \geq 0$ for all values $2 \leq s \leq 17$ as required; and only for $s = 4$ does equality obtain.

| s | 2 | 3 | 4 | 5 | 6 | 7 | 8 | 9 | 10 | 11 | 12 | 13 | 14 | 15 | 16 | 17 |
|---|---|---|---|---|---|---|---|---|----|----|----|----|----|----|----|----|
| $j_0$ | 1 | 1 | 2 | 2 | 2 | 3 | 3 | 3 | 4 | 4 | 4 | 5 | 5 | 5 | 6 | 6 |
| 5s/2 | 5 | 7.5 | 10 | 12.5 | 15 | 17.5 | 20 | 22.5 | 25 | 27.5 | 30 | 32.5 | 35 | 37.5 | 40 | 42.5 |
| 7j | 7 | 7 | 14 | 14 | 14 | 21 | 21 | 21 | 28 | 28 | 28 | 35 | 35 | 35 | 42 | 42 |
| d | 2 | 4.5 | 0 | 2.5 | 3 | 0.5 | 3 | 4.5 | 1 | 3.5 | 6 | 1.5 | 4 | 6.5 | 2 | 4.5 |

#### 2.1.2.3  $s = 1$

When $s = 1$, the picture is as follows (an example is $A = \{1, 31, 47\}$ with $n = 15$, $p = 2$):

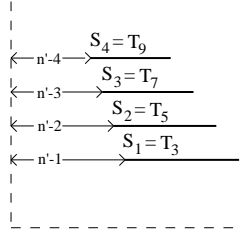

We have the following exact formulae:

$str(S_j) = n' - j$
$len(S_j) = n' - j$
$ord(S_j) = 2j + 1$

$S_j$ is certainly covered by $T_0$ when $str(S_j) + len(S_j) \leq n'+1 \iff 2n' - 2j \leq n' + 1 \iff 2j \geq n'-1$; we now consider two cases according as $n'$ is even or odd.

a) $n'$ even, $n' = 2m$

Choose $j_0 = m \Rightarrow 2j_0 = 2m \geq n'-1 = 2m-1$; so $S_{j_0}$ is covered by $T_0$.
$a_2 = 2n + 1 = 4m + 1$; so $ord(S_{j_0}) = 2m + 1 \leq a_2/2 + 1$ as required.

b) $n'$ odd, $n' = 2m+1$

Choose $j_0 = m \Rightarrow 2j_0 = 2m \geq n'-1 = 2m$; so $S_{j_0}$ is covered by $T_0$.
$a_2 = 2n + 1 = 4m + 3$; so $ord(S_{j_0}) = 2m + 1 < a_2/2$ as required.

### 2.1.3 $s > n'/2$

In this case, the threads with offset < n' group as follows:

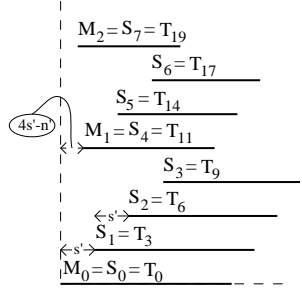

We write $s' = n'-s \Rightarrow s = n'-s'$, $a_2 = 3n'-s'$. Note that we do not include $T_1$ (since this has offset $= n'$); the groups are $\{T_0, T_3, T_6, T_9\}$, $\{T_{11}, T_{14}, T_{17}\}$, $\{T_{19}, T_{22}, ...\}$ etc.

The number of threads in each group is determined by the offset of the group's minimal thread, which satisfies $0 \leq$ offset $< s'$; since $n' = ks + t$, there are either $k$ or $k+1$ threads in each group.

The difference in order between the threads in each group is 3, and the difference in order between the last thread of one group and the first thread of the next group is 2; so the difference in order between successive minimal threads is either $3k-1$ or $3k+2$.

The difference in length between successive threads in the same group is 2, and between the last thread of one group and the first thread of the next is 1; so the difference in length between successive minimal threads is either $2k-1$ or $2k+1$.

The first group always contains $k+1$ threads (since $str(S_0) = 0$), so we have the following bounds:

$ord(M_j) \leq j(3k+2)$ - (32)
$len(M_j) \leq (n'+1) - (2k+1) - (j-1)(2k-1) = (n'-2k) - (j-1)(2k-1)$ - (33)
$str(M_j) < s'$ - (34)

We write $n' = ks' + t$, $0 \leq t < s'$, and now show that we may assume $t > 0$; for if $t = 0$, $str(M_1) = 0$ and so $M_1$ is covered by $T_0$; we have:

$ord(M_1) = 3k+2$, and so $q \leq 3k+1$
$a_2 = 3n' - s' = 3ks' + 3t - s' \geq (3k-1)s'$

So for $s' \geq 3$, $a_2 - 2q \geq (9k-3) - 2(3k+1) = 3k-5 > 0$, since $s > n'/2 \Rightarrow s' < n'/2 \Rightarrow k \geq 2$.
So $a_2 - 2q > 0 \Rightarrow q < a_2/2$ as required - provided that $s' \geq 3$.
We deal with the cases $s' = 1$, $s' = 2$ specially below, so we may otherwise assume that $t > 0$.

For $s' \geq 3$, we show $q \leq a_2/2$ by finding a minimal thread $M_{j_0}$ that is covered by $T_0$ and which satisfies $ord(M_{j_0}) \leq a_2/2 + 1$. From (33) and (34), $M_j$ is certainly covered by $T_0$ when*:

$len(M_j) + s' \leq n'+1 \iff (n'-2k) - (j-1)(2k-1) + s' \leq n'+1 \iff s' \leq (j-1)(2k-1) + (2k+1) = j(2k-1) + 2 \iff j \geq (s'-2)/(2k-1)$ - (35)

We write $d = a_2/2 + 1 - ord(M_j)$:

$d \geq 3n'/2 - s'/2 + 1 - j(3k+2) = 3(ks' + t)/2 - s'/2 + 1 - j(3k+2)$
So: $2d \geq 3(ks' + 1) - s' + 2 - 2j(3k+2) = s'(3k-1) + 5 - 2j(3k+2)$

Now (35) is satisfied by some integer $j < (s'-2)/(2k-1) + 1$, so:

$2(2k-1)d \geq s'(2k-1)(3k-1) + 5(2k-1) - 2(s'-2)(3k+2) - 2(2k-1)(3k+2)$
$= s'(6k^2 - 5k + 1) + (10k-5) - 2s'(3k+2) + 4(3k+2) - 2(6k^2 + k - 2)$
$= s'(6k^2 - 11k - 3) - (12k^2 - 20k - 7)$ - (36)

[ * Note that $len(M_j) + (s'-1) \leq n'+1$ is sufficient, but this does not avoid the $k = 2$ issue. ]

### 2.1.3.1 $s' \geq 3$

Infuriatingly, in (36) we find $(6k^2 - 11k - 3) = -1$ for $k = 2$; so we must treat $k = 2$ as a a special case.

### 2.1.3.1.1 $k > 2$

For $k \geq 3$, $(6k^2 - 11k - 3) > 0$, and so for $s' \geq 3$ we have from (36):
$2(2k-1)d \geq 6k^2 - 13k - 2 \geq 0$ for all $k \geq 3$.
So in this case $ord(M_{j_0}) \leq a_2/2 + 1 \Rightarrow q \leq a_2/2$ as required.

*2.1.3.1.2  k = 2*
When k = 2, the threads group as follows (cf diagram in 2.1.3):

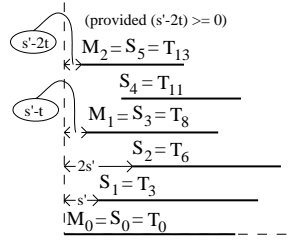

Each group contains either 2 or three threads, eg $\{S_0, S_1, S_2\}$, $\{S_3, S_4\}$. We have:
$\quad$ str($M_1$) = 3s'-n' = 3s' - (2s'+t) = s'-t
$\quad$ str($M_2$) = 5s'-2n' = 5s' - (4s'+2t) = s'-2t
$\quad$ ....
So in general we have:
$\quad$ str($M_j$) = (-jt) mod s'
Each time s' is added, there is a 3-thread group; so if we write:
$\quad$ fs' - jt = g $\quad\quad$ where $0 \le g < s'$
then there are exactly f 3-thread groups before $M_j$, and so we have (cf (32)):
$\quad$ ord($M_j$) = f(3k+2) + (j-f)(3k-1) = 8f + 5(j-f) = 3f + 5j
From (35) we know that $M_j$ is certainly covered by $T_0$ when $j \ge (s'-2)/3$; we now write s' = 3u + x for some $0 \le x \le 2$, and choose $j_0 = u$: we see that $j_0 \ge (s'-2)/3$ in all cases, and so $M_{j0}$ is certainly covered by $T_0$. We have:
$\quad$ fs' - $j_0$t = g  => f(3u+x) - ut = g  => f = g/(3u+x) + ut/(3u+x) < 1 + t/3  since $0 \le g < s' = 3u+x$
So:$\quad$ ord($M_{j0}$) = 3f + 5$j_0$ < t + 3 + 5u
and:$\quad$ $a_2$ = 5s' + 3t = 5(3u+x) + 3t $\ge$ 15u + 3t
So $a_2$ - 2*ord($M_{j0}$) > (15u + 3t) - 2(5u + t + 3) = 5u + t - 6 $\ge$ 0  since t > 0, and s' $\ge$ 3  => u $\ge$ 1.
So ord($M_{j0}$) $\le a_2/2$ => q < $a_2/2$ as required.

*2.1.3.2  s' = 2*
In this case, $S_1 = T_3$ is covered by $T_0$, since str($T_3$) = s' = 2, and len($T_3$) = (n'+1) - 2. So q $\le$ 2, and $a_2$ = 3n' - s' = 3(ks' + t) - s' $\ge$ 6k - 2 $\ge$ 10; so q < $a_2/2$ as required.

*2.1.3.3  s' = 1*
As above, $S_1 = T_3$ is covered by $T_0$, so q $\le$ 2. $a_2$ = 3(ks' + t) - s' $\ge$ 3k - 1 $\ge$ 5; so q < $a_2/2$ as required.

*2.2  p $\ge$ 3*
We have already shown that n < $a_2/2$, and experiment suggests that $2q_{max} \le a_2 + 1$ for all non-canonical stride generators meeting these conditions (that is, with $C_2 = 1$, $C_1 > a_2/2$ and p $\ge$ 3) - but not for canonical ones, where $2q_{max} = a_2 + 2$ is observed. (Here $q_{max}$ is the order of the smallest possible thread in the stride generator, this guaranteeing that q $\le q_{max}$). In detail, we find:
$\quad$ A = {1, 30, 58} is an example of a canonical fundamental stride generator with $2q_{max} = a_2 + 2$; $q_{max}$ = 16, n = 15, p =14, n' = 2.
$\quad$ Examples of non-canonical stride generators with $2q_{max} \ge a_2$ are:
$\quad\quad$ n' = 10 $\quad$ {1, 31, 52} $\quad\quad$ $q_{max}$ = 16 $\quad\quad$ n = 11,  p = 3
$\quad\quad$ n' = 10 $\quad$ {1, 32, 54} $\quad\quad$ $q_{max}$ = 16 $\quad\quad$ n = 11,  p = 3
$\quad\quad$ n' = 11 $\quad$ {1, 34, 57} $\quad\quad$ $q_{max}$ = 17 $\quad\quad$ n = 12,  p = 3
$\quad\quad$ n' = 11 $\quad$ {1, 35, 59} $\quad\quad$ $q_{max}$ = 17 $\quad\quad$ n = 12,  p = 3
[ This pattern repeats. In fact, $q_{max}$ for p $\ge$ 3 is at a maximum when p = 3 and n' ~ $a_2/3$ or n' = 2, and dips in between these two values. ]
We are able to show that $q_{max} < a_2/2 + 1$ for all non-canonical stride generators for p $\ge$ 4, but not for p = 3: although it seems to be true, we can only manage q < $a_2/2$ + (5/2)!

Using the "thread length" argument from section 2.4.2 above, we have:
$\quad$ len($T_i$) $\le$ len($T_0$) - i$C_2$ + (i-1)/p = n + 1 - i + (i-1)/p
So len($T_i$) $\le$ 0 when:
$\quad$ n + 1 - i(1 - 1/p) - 1/p $\le$ 0 <=> i(1 - 1/p) $\ge$ n + (1 - 1/p) <=> i $\ge$ pn/(p-1) + 1 => q < np/(p-1) + 1 $\quad\quad$ - (37)
Substituting n = n' + p - 2 and using pn' < $a_2$ (see 2.4.3), we have:
$\quad$ q < (n' + p - 2)p/(p-1) + 1 < ($a_2$/p + p - 2)p/(p-1) + 1
So q < $a_2/2$ + 1 when ($a_2$/p + p - 2)p/(p-1) $\le a_2/2$ <=> 2p($a_2$/p + p - 2) $\le a_2$(p-1) <=> 2$a_2$ + 2$p^2$ - 4p - $a_2$p + $a_2$ $\le$ 0
$\quad$ <=> $\quad\quad$ 2$p^2$ - ($a_2$ + 4)p + 3$a_2$ $\le$ 0 $\quad\quad$ - (38)
So to complete the proof, we have only to show that (38) is true.

From (4) in 2.4.1 we have:  1 < $a_2$/(p(p+1)) + 1/(p+1) <=> p(p+1) < $a_2$ + p <=> $a_2 > p^2$ $\quad\quad$ - (39)
$\quad\quad$ [which is an interesting result in its own right: but remember that this is true for fundamental non-canonical stride generators only]

*2.2.1  p $\ge$ 5*
Using (39) we obtain:
$\quad$ 2$p^2$ - ($a_2$ + 4)p + 3$a_2$ < 5$a_2$ - ($a_2$ + 4)p = (5 - p)$a_2$ - 4p < (5 - p)$a_2$ $\le$ 0 for all p $\ge$ 5
So (38) is true when p $\ge$ 5, as required.

*2.2.2  p = 4*
When p = 4, n' < $a_2$/4 and so n < $a_2$/4 + 2. Substituting directly in (37) gives:
$\quad$ q < 4($a_2$/4 + 2)/3 + 1 = ($a_2$ + 11)/3
So n + q < ($a_2$ + 8)/4 + ($a_2$ + 11)/3 = (7$a_2$ + 68)/12  which is $\le a_2$ when 5$a_2 \ge$ 68 => $a_2 \ge$ 14; but (39) requires $a_2$ > 16, so this case is proved.

An alternative approach allows us to improve on this by actually showing that q < $a_2/2$ + 1; we substitute directly in (38):
$\quad$ 2$p^2$ - ($a_2$ + 4)p + 3$a_2$ = 32 - 4($a_2$ + 4) + 3$a_2$ = 16 - $a_2$ < 0 since (39) requires $a_2$ > 16.

*2.2.3  p = 3*

*Page 45*

When p = 3, n' < $a_2/3$ and so n < $a_2/3$ + 1. Substituting directly in (37) gives:
   q < 3($a_2$/3 + 1)/2 + 1 = ($a_2$ + 5)/2
So n + q < ($a_2$ + 3)/3 + ($a_2$ + 5)/2 = (5$a_2$ + 21)/6 which is ≤ $a_2$ when $a_2$ ≥ 21; (39) requires $a_2$ > 9, so we have only to consider 10 ≤ $a_2$ ≤ 20; the results are as follows (the "closest" we get is $q_{max}$ = $a_2$/2 + 1 for {1, 19, 32}):

| $a_2$ | 10 | 11 | 12 | 13 | 14 | 15 | 16 | 17 | 18 | 19 | | 20 |
|---|---|---|---|---|---|---|---|---|---|---|---|---|
| $a_2$/4 | 2.5 | 2.75 | 3 | 3.25 | 3.5 | 3.75 | 4 | 4.25 | 4.5 | 4.75 | | 5 |
| $a_2$/3 | 3.33 | 3.66 | 4 | 4.33 | 4.66 | 5 | 5.33 | 5.66 | 6 | 6.33 | | 6.66 |
| n' (*) | 3 | 3 | - | 4 | 4 | 4 | 5 | 5 | 5 | 5 | 6 | 6 |
| $a_3$ | 17 | 19 | - | 22 | 24 | 26 | 27 | 29 | 31 | 33 | 32 | 34 |
| n | 4 | 4 | - | 5 | 5 | 5 | 6 | 6 | 6 | 6 | 7 | 7 |
| $q_{max}$ ** | (C) | (C) | | (C) | (C) | (C) | 8 | (C) | (C) | (C) | 10 | 10 |

  * This line gives possible values for n'; remember that $a_2$/4 < n' < $a_2$/3
 ** (C) indicates that the stride generator is canonical